\documentclass[11pt]{amsart}

\usepackage[english]{babel}
\usepackage[T1]{fontenc}
\usepackage[utf8]{inputenc}
\usepackage{amsbsy}
\usepackage{amsmath}
\usepackage{amssymb}
\usepackage{amsfonts}
\usepackage{amsthm}
\usepackage{bm}
\usepackage{graphicx}
\usepackage{hyperref}
\usepackage[active]{srcltx}
\usepackage{upgreek}
\usepackage{xcolor}
\usepackage{mathrsfs}

\usepackage{mathtools}

\usepackage{enumitem}

\usepackage{tikz}

\usepackage{imakeidx}
\makeindex
 
\newcommand{\m}[1]{
\ifdefequal{#1}{1}
{\mathbbm{#1}}
{\mathbb{#1}}
}
\newcommand{\q}[1]{\mathscr{#1}}

\newcommand{\N}{\mathbb{N}}

\newcommand{\R}{\mathbb{R}}
\renewcommand{\S}{\mathbb{S}}

\newcommand{\Z}{\mathbb{Z}}

\renewcommand{\ge}{\geqslant}

\newcommand{\gm}{\mathfrak{m}}

\newcommand{\grad}{\nabla}
\newcommand{\Rdeux}{{\mathbb{R}^2}}

\newcommand{\heff}{h_{\text{eff}}}

\newcommand{\ii}{\infty}
\newcommand{\dom}{\mathcal{D}om}

\newcommand{\loc}{\mathrm{loc}}

\newcommand{\norm}[2]{\left\Vert #1\right\Vert_{#2}}

\newcommand{\normLii}[1]{\left\Vert #1\right\Vert_{L^\ii}}
\newcommand{\normLdeux}[1]{\left\Vert #1\right\Vert_{L^2}}

\newcommand{\normHun}[1]{\left\Vert #1\right\Vert_{H^1}}
\newcommand{\psLdeux}[2]{\left\langle #1,#2 \right\rangle_{L^2}}
\newcommand{\psLdeuxLdeux}[2]{\left\langle #1,#2 \right\rangle_{L^2\times L^2}}
\newcommand{\vvvert}[1]{\big\vert\kern-0.25ex\big\vert\kern-0.25ex\big\vert #1\big\vert\kern-0.25ex\big\vert\kern-0.25ex\big\vert}

\newcommand{\spanned}[1]{\mathrm{Span}(#1)}
\newcommand{\normLdeuxLdeux}[1]{\left\Vert #1\right\Vert_{L^2\times L^2}}

\newcommand{\normR}[1]{\left| #1\right|}

\newtheorem{thm}{Theorem}[section]

\newtheorem{claim}[thm]{Claim}

\newtheorem{lem}[thm]{Lemma}
\newtheorem{prop}[thm]{Proposition}

\newtheorem*{thm*}{Theorem}

\newtheorem{rem}[thm]{Remark}

\newtheorem*{notations}{Notations}

\theoremstyle{remark}

\title[Asymptotic stability of skyrmions for (LLG)]{Asymptotic stability of a Bloch skyrmion for the Landau-Lifshitz-Gilbert equation}

\date{}

\author[Jordan Berthoumieu]{Jordan Berthoumieu}
\address{Karlsruhe Institute of Technology (KIT), Department of Mathematics, Englerstraße 2, 76131 Karlsruhe Germany}
\email{jordan.berthoumieu@kit.edu}

\author[Raphaël Côte]{Raphaël Côte}
\address{IRMA UMR 7501, Université de Strasbourg, CNRS, F-67000 Strasbourg, France \newline
\& USIAS, Université de Strasbourg, F-67083 Strasbourg France}
\email{rcote@unistra.fr}

\allowdisplaybreaks

\begin{document}

\begin{abstract}
We prove the asymptotic stability of isolated skyrmions under the flow of the two-dimensional Landau-Lifshitz-Gilbert equation. We consider a micromagnetic energy functional involving the exchange energy, the Dzyaloshinskii-Moriya interaction, and the Zeeman interaction. We are interested in the stability of chiral Bloch skyrmions, which are local minimizers of this energy. We prove that data in $H^2$ which are near in the skyrmion in a slighlty stronger  than $H^1$ sense, converge with exponential rate, in $H^2$, to the skyrmion, up to a change of gauge.
\end{abstract}

\maketitle

\tableofcontents


\numberwithin{equation}{section}

\section{Introduction}

\subsection{Micro-magnetic energy and skyrmions}

Skyrmions are topologically non-trivial field configurations, appearing in physics literature as spin textures. They were first introduced by Tony Skyrme in 1962~\cite{Skyrme1}, and since then, several mathematical models have been proposed, depending on its intrinsic properties and physical mechanism under consideration on the chosen mathematical framework.  They appear in various fields of very different nature at nanoscale, including Bose-Einstein condensates~\cite{AlKhStoo}, or nematic liquid crystals~\cite{BogdaRossShes}.  These nano-structures have been experimentally produced during the past two decades~\cite{VigoGui}, and stand a fair chance of success in their potential industrial and technological applications~\cite{BogKis}.

Skyrmions may be viewed as topological solitons~\cite{KosIvKo1}, since they are localized stationary solutions whose profiles are protected by a topological invariant called the \textit{skyrmion number}, see~\eqref{skyrmion number}\index{skyrmion number} below; and they seem to be energetically favored structures~\cite{BogdaHub2}. More specifically, given two domains  $\mathcal{D}_1, \mathcal D_2$ where $\mathcal D_1 \in\{\S^2,\S^3,\Omega\}$, with $\Omega$ an open subset (not necessarily bounded) of $\R^2$ or $\R^3$, and $\mathcal{D}_2\in\{\R^3,\S^2,\S^3\}$, skyrmions can be defined as energy minimizing maps $m:\mathcal{D}_1\rightarrow\mathcal{D}_2$ constrained to a fixed non-trivial topological degree, whenever the source and target domains have the same dimension.

In this article, we give special interest to skyrmions in \textit{micromagnetics}, as stable magnetic textures in ultrathin magnetic films or layers~\cite{RomHanneMenetal,NagaNaoto}. The magnetic field in such ferromagnetic material shall be described by a map defined on the whole plane, with values in $\m S^2$:
\[ m=(m_1,m_2,m_3):\R^2 \rightarrow \S^2 \subset \m R^3. \]
We denote $(e_1,e_2,e_3)$ the canonical basis of $\m R^3$, and we focus on the following micromagnetic energy functional
\begin{align*}\label{expression énergie}
    E_{\kappa,h}(m) &:=E_{ex}(m) + \kappa E_{DMI}(m) + hE_Z(m) \\
    & :=\int_{\R^2}\bigg(\dfrac{|\grad m|^2}{2}+\kappa m \cdot \grad\times m + \dfrac{h}{2}|m-e_3|^2\bigg) dx.
\end{align*}
In this framework, $\kappa,h>0$ are positive constants. The gradient $\grad$ acts according to the two-dimensional setting and should be viewed as $\grad=(\partial_1,\partial_2,0)$. In particular, we have $\grad\times m =(\partial_2 m_3,-\partial_1 m_3,\partial_1 m_2 - \partial_2 m_1) $ and $\cdot$ designates the usual euclidian scalar product\index{euclidian scalar product $\cdot$ on $\R^3$} on $\R^3$. $E_{ex}$ is the Dirichlet (exchange) energy and models the interaction with nearest neighbors. $E_{DMI}$ is the (bulk) Dzyaloshinskii-Moriya interaction that accounts for the spin-orbit coupling effect, with $\kappa$ being the coupling constant. And $E_Z$ is the Zeeman energy, characterizing the effect of an external oriented (applied) magnetic field which tends to align the spins in the direction $e_3$, with constant amplitude $h$.

Notice that $\m R e_3$ is the easy axis, and so magnetic textures with finite energy  may be viewed as excited states of the trivial background state $m\equiv e_3$, which is the least energy configuration for $E_{\kappa,h}$. 

In the literature, the micromagnetic energy can take various forms and involve also other terms than in $E_{\kappa,h}$. For example, an extra anisotropy term, or a stray field term~\cite{BogdaYablon1,BernMuraSim1}. We will focus on $E_{\kappa,h}$ here, as it already concentrates the main features and difficulties. 

\bigskip

In the sequel, we shall show that the energy $E_{\kappa,h}$ is well-defined on the natural energy set of Ginzburg-Landau type given, for $s\geq 1$, by
\begin{equation*}
    H^s_* : = \Big\{ m\in H^1_{\loc}(\R^2,\R^3): |m|=1 \ \ a.e. \text{ and }\  m - e_3 \in H^s(\R^2,\R^3)\Big\}.
\end{equation*}

This set is not a subspace, in particular due to the (non-convex) constraint to be $\m S^2$-valued. However, it is naturally a (complete) metric space when endowed with the natural metric $\norm{m-\tilde m}{H^s}$.
Moreover, the formal identity
\begin{equation*}
    \int_\Rdeux m \cdot(\grad\times m) dx= \int_\Rdeux (m-e_3) \cdot(\grad\times m) dx,
\end{equation*}
can be made rigorous for $m \in H^1_*$ (c.f.\cite{Melcher1}), and this allows to define properly the functional $E_{\kappa,h}$ on $H^1_*$. 

\bigskip

On a theoretical level, skyrmions were first rigorously shown to exist using variational methods: they were defined as critical points of the energy $E_{\kappa,h}$ in a non-trivial homotopy class characterized by the quantity 
\begin{equation}\label{skyrmion number}
    Q(m)=\dfrac{1}{4\pi}\int_{\R^2} m.(\partial_1 m\times \partial_2 m).
\end{equation}
This quantity is well defined and takes values in $\m Z$ for $m \in H^1_*$; it is called the \emph{skyrmion number}\index{skyrmion number}.

The first mathematical evidence of the existence of skyrmions is given in~\cite{Melcher1}, where the author proved the existence of minimizers for the energy with given skyrmion number $Q=-1$ and for a special range of parameters $\kappa$ and $h$. Let us underline the role of the Dzyaloshinskii-Moriya interaction in this matter: due to its symmetry-breaking features (see \cite{BogdaYablon1}), this term is responsible for the emergence of \textit{magnetic chiral skyrmions} in ferromagnetic heterostructures~\cite{RomHanneMenetal,BogKudYablon,BogHub1}. 

We refer to \cite{MuraSimSlas1} for the existence of skyrmions with any prescribed negative skyrmion number on a bounded domain, and to \cite{BernMuraSim1,GustafWang,DorMel,IbraShimi} for works which provide a quantitative analysis to refine the variational characterization of the skyrmion and their relation with harmonic maps.

\subsection{Chiral axisymmetric skyrmions}

Henceforth, we set $\kappa=1$, and we focus on the corresponding $E:=E_{1,h}$, which fits the framework of \cite{LiMelcher}. $E_{ex}$ and $E_{Z}$ are invariant by any rotation or reflexion. However, this is not the case of $E_{DMI}$. Nevertheless the symmetry group is not trivial and we are led to consider \emph{chiral axysymmetric} profile of degree $N \in \m Z$, i.e. magnetizations $m$ such that for any rotational matrix $R\in SO(2)$, and any $x\in\R^2$, 
\begin{equation}\label{SO(2) equivariance}
    m(Rx)=\widetilde{R}^{-N} m(x)=\widetilde{^t R}^N m(x) \quad \text{with } \widetilde{R}=\begin{pmatrix}
    R & 0\\
    0 & 1
\end{pmatrix}, \index{$\widetilde{R}$}
\end{equation}
and
\begin{equation}\label{combined reflection}
    m(Fx)=\widetilde{-F}m(x) \quad \text{with } F=\begin{pmatrix}
    1 & 0\\
    0 & -1
\end{pmatrix}, \index{$F$}
\end{equation}

We refer to Appendix~\ref{appendix: symétries} and to Subsection 1.3 in~\cite{GustafWang} for more details on these symmetries\footnote{One can also find the terminology "co-rotational" for maps satisfying both previous symmetries. In other words, these symmetries correspond to the combination of the $SO(2)$-equivariance and~\eqref{combined reflection}.}. We represent any $x\in\R^2$ with its polar coordinates $(r,\psi)\in [0,+\ii)\times [0,2\pi)$. General axysymmetrical chiral skyrmions take the form
\begin{equation}\label{définition skyrmion_0}
    m^{(N,\gamma)}_*(r,\psi)=\begin{pmatrix}
        \sin\big(\theta(r)\big)\cos(\gamma -N\psi)\\
        \sin\big(\theta(r)\big)\sin(\gamma -N\psi)\\
        \cos\big(\theta(r)\big)
    \end{pmatrix} ,
\end{equation}

\begin{figure}[h]
\centering
\includegraphics[scale=0.5]{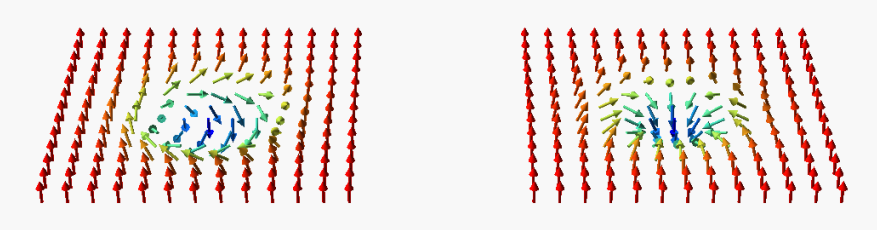}
\caption{\textit{Schematic sketch of a Bloch (left) and a Néel skyrmion, with skyrmion number $Q=-1$ (courtesy of the authors~\cite{SchulYuXiao}).}}
\end{figure}

They are therefore completely determined by the polar-angle $\theta: \m R_+ \to \m R$ given in Proposition~\ref{proposition: propriétés du skyrmion} below; the constant $\gamma$ is called the \emph{helicity}\index{helicity} of the skyrmion, and $N=Q(m_*)$ is the \emph{skyrmion number}\index{skyrmion number}. It is not known whether the minimizers of $E$ with a given skyrmion number obtained in~\cite{Melcher1,MuraSimSlas1} are of the above form~\eqref{définition skyrmion_0}, and this goes beyond the scope of this article. However, it is proven in~\cite{LiMelcher} that within the equivariant class \eqref{définition skyrmion_0}, any non trivial critical point $m_*^{(N,\gamma)}$ of the micromagnetic energy $E$ necessarily satisfies $N=-1$, and that $E_{DMI}(m_*)$ is minimized for $\gamma=\frac{\pi}{2}$. This particular solution is called the \emph{Bloch skyrmion}\index{Bloch/Néel skyrmion}\footnote{The skyrmion with $N=-1$ and $\gamma=0$ is called the Néel skyrmion, or hedgehog skyrmion and we refer to Appendix~\ref{remarque: Néel skyrmion} for more details.} and henceforth we write
\begin{equation}\label{définition skyrmion}
    m_*(r,\psi):=m_*^{(-1,\frac{\pi}{2})}(r,\psi)=\begin{pmatrix}
        -\sin\theta(r)\sin\psi\\
        \sin\theta(r)\cos\psi\\
        \cos\theta(r)
    \end{pmatrix} = \sin\theta(r) e_\psi + \cos\theta(r) e_3,
\end{equation}
with $(e_r,e_\psi,e_3)$ is the cylindrical coordinate system, given by
\begin{equation*}
    e_r=\begin{pmatrix}
        \cos\psi \\
        \sin\psi\\
        0
    \end{pmatrix}\quad \text{and}\quad e_\psi = \begin{pmatrix}
        -\sin\psi\\
        \cos\psi\\
        0
    \end{pmatrix}.
\end{equation*}

As we shall see it in the next subsection, static solutions such as $m_*$ are shown to be linearly stable critical points of $E$: indeed they solve 
\begin{equation}\label{euler lagrange equation point critique}
    m_*\times \delta E(m_*)=0,
\end{equation}
where $\delta E(m_*)$ designates the gradient given in the variational sense for any $m\in H^2(\R^2)$ by the formula
\begin{equation*}
    \delta E(m)=-\Delta m +2\grad\times m + h(m-e_3).
\end{equation*}

Equation~\eqref{euler lagrange equation point critique} reduces to an ODE on the profile $\theta$ of the skyrmion, c.f.~\cite{LiMelcher}. The main properties satisfied by the profile of the skyrmion are summarized in the following proposition.
\begin{prop}[\cite{LiMelcher}]\label{proposition: propriétés du skyrmion}
    There exists $h_0>1$ large enough such that for any $h\geq h_0$, the ODE 
\begin{equation}\label{équation sur theta}
    \theta'' + \dfrac{\theta'}{r} - \dfrac{\sin\theta\cos\theta}{r^2}+ \dfrac{2\sin^2(\theta)}{r}-h\sin\theta=0,
\end{equation}
with boundary conditions $\theta(0)=\pi$ and $\theta(r)\rightarrow 0$ as $r\rightarrow +\ii$ has a unique solution $\theta$. Moreover, $\theta\in \q C^\ii\big(\R^*_+,(0,\pi)\big)$ is monotonically decreasing. Finally, we have the asymptotics
\begin{align*}
    \theta(r) & = \left( \dfrac{\alpha_h}{\sqrt{r}} + o(1) \right)e^{-\sqrt{h}r}  \quad \text{as} \quad r \to +\infty, \quad\text{and} \\
    \theta(r) & = \pi - \dfrac{h}{2}r + o(r) \quad \text{as} \quad r \to 0,
\end{align*}
for a positive constant $\alpha_h$. In particular $m_*\in H_*^\ii=\bigcap_{s\geq 0}H^s_*$. 
\end{prop}

In the sequel, we assume that $h$ satisfies the hypothesis of Proposition~\ref{proposition: propriétés du skyrmion}. We then turn to the dynamical analysis of the skyrmion $m_*$. Note that $E$ is invariant by translation, thus $m_*\big(.-y\big)$ is also a critical point of the energy for any $y\in \R^2$. This symmetry has important consequences from the spectral point of view.

\subsection{The Landau-Lischitz-Gilbert flow}

We are interested in the dynamics driven by the Landau-Lifshitz-Gilbert equation
\begin{equation}\label{llg}\tag{LLG}
    \partial_t m =m\times \heff(m) -\alpha m \times \big(m\times \heff(m)\big),
\end{equation}
with an effective magnetic field given by the opposite of the gradient of the energy, namely
\begin{equation}\label{expression heff}
    \heff(m):=-\delta E(m)=\Delta m-2\grad\times m -h(m-e_3) .
\end{equation}
This evolution equation was originally introduced in~\cite{LandLif1} to model magnetization dynamics in ferromagnetic materials. In the latter article, both cases $\alpha=0$ and $\alpha>0$ are derived, depending on the kind of interaction they consider between the moments in the ferromagnetic crystal. To consider a positive relaxation constant $\alpha >0$ amounts to taking into account the damping effect of spin waves~\cite{GuoDing0}. 

The first issue to address regarding~\eqref{llg} is local well-posedness. One can prove a result for general finite energy data, especially as in space dimension at least 2, $H^1$ does not provide a pointwise bound (relevant for the constraint to be $\m S^2$-valued). Nevertheless, one is often interested in the dynamics near a specific solution. In order to provide a unified framework (with respect to other results such as~\cite{CarbLabb1,CotIgnat} for domain walls in nanowires or in more general frameworks~\cite{shimizu2026,GuoDing0}) we consider data localized around a fixed back-ground $\mathfrak m$, suitable in a sense made precise below.

For $s\geq 1$, an appropriate set for the background $\mathfrak m$ is the so-called Zhidkov space
\begin{equation*}
    \mathcal{Z}^s(\R^d):=\big\{ \mathfrak m\in L^\ii(\R^d)\big| \grad \mathfrak m\in H^{s-1}(\R^d)\big\}. \index{$\mathcal{Z}^{s}(\R^d)$}
\end{equation*}
Notice that for $d=2$ or $3$ and $s > 1$,  $\mathcal{Z}^s(\R^d)$ and $H^s(\m R^d)$ embed continuously in $\q C^0(\m R^d)$. All our statements will be made at that level of regularity, so that the target space be $\m S^2$ can be understood in the pointwise sense both in time and space, and this will not raise any issue.

Here is our first statement.

\begin{thm}\label{théorème: local well-posedness}
Let $s\geq 2,d\in\{2,3\}$ and $\mathfrak{m}\in \mathcal{Z}^{s}(\R^d)$ satisfying
\begin{equation}\label{Z_0 in H^s-1}
    \mathfrak{m}\times \heff(\mathfrak{m}) \in H^{s-1}(\R^d).
\end{equation}
Then, the equation \eqref{llg} is well-posed on $\mathfrak{m} + H^s(\R^d)$, as precised below.

Let $w_0\in H^s(\R^d)$. There exist a maximal time of existence $T_{\max}>0$ and a unique function 
\[ w\in \q C^0\big([0,T_{\max}),H^s(\R^d)\big)\cap L^2\big([0,T_{\max}),H^{s+1}(\R^d)\big) \cap \q C^1([0,T_{\max}),H^{s-2}(\R^d)\big) \]
 such that $m :=\mathfrak{m} + w$ solves~\eqref{llg} with 
\[ m(0) := \mathfrak{m} + w_0. \]
If we have initially $|m(0,x)|=1$ for all $x \in \m R^d$, then for any $t\in [0,T_{\max})$ and $x \in \m R^d$, $|m(t,x)|=1$. Furthermore, if $T_{\max}<\ii$, then 
\begin{equation*}
    \norm{w(t)}{H^s} \to +\ii \quad \text{as} \quad t \to T_{\max}.
\end{equation*}
Finally, the flow map $w_0\mapsto w$ is continuous on $H^{s-1}(\R^d)$: for $T < T_{\max}$, the solution with data $\tilde w_0$ in a neighbourhood $\q V$ of $w_0$ in $H^{s-1}(\m R^d)$ has a lifespan at least $T$, and the map $\tilde w_0 \mapsto w$ is continuous from $\q V$ to $\q C^0([0,T],H^{s-1}(\R^d)\big)$.
\end{thm}


As a consequence, we infer a local well-posedness for localized perturbation of the skyrmion $m_*$, for which condition \eqref{Z_0 in H^s-1} is clearly satisfied (the quantity evaluates to $0$, c.f. \eqref{euler lagrange equation point critique}). It also provides the energy dissipation identity, allows to justify the manipulation on time-dependent quantities that we use in Subsections~\ref{section: dissipation term} and~\ref{section: higher order term}.

\begin{thm}\label{théorème: LWP of the Cauchy problem}
Let $s\in\{2,3\}$ and $m_0\in H^s_*$. There exists a maximal time $T_{\max}>0$ and a solution $m$ to~\eqref{llg} such that $m-m_*\in \q C^0\big([0,T_{\max}), H^s(\R^2)\big)$ with initial data $m(0)=m_0$. Moreover $t\mapsto E\big(m(t)\big)\in \q C^1\big([0,T_{\max}),\R\big)$ and we have the following dissipation identity
\begin{equation*}
    \dfrac{d}{dt}\Big[E\big(m(t)\big)\Big]=-\alpha\int_{\R^2} \Big|m(t)\times \delta E\big(m(t)\big)\Big|^2.
\end{equation*}
Finally, if $m_0\in H^3_*$, $t\mapsto \normLdeux{\delta E\big(m(t)\big)}^2\in \q C^1\big([0,T_{\max}),\R\big)$.
\end{thm}

The proofs of both Theorems \ref{théorème: local well-posedness} and \ref{théorème: LWP of the Cauchy problem} are done in  Appendix~\ref{appendix: LWP}.

\subsection{Main result: asymptotic stability of skyrmions}

The stability question has been handled, at the stationary level, in~\cite{LiMelcher} and the analysis in this article relies substantially on this result. As a matter of fact, the coercive properties showed there provide the structure needed to prove asymptotic stability. Namely, we prove that any initial data $m_0$ taken initially close to a skyrmion $m_*$ yields a solution $m(t)$ to~\eqref{llg}, globally defined in time, which eventually converges to the same skyrmion as time tends to $+\ii$, up to a time-dependent translation. We also mention the work~\cite{IbraShimi2}, where stationary linear stability is proven in a more general framework, where a wide range of parameters (exchange, DMI, Zeeman and easy-plane anisotropy) is considered. Note that in both aforementioned articles, a more geometrical point of view is adopted, since perturbations with values in the tangent plane of $m_*$ are considered. This differs strongly from our approach, where we consider perturbations valued in $\R^3$, including the component in the normal direction of the sphere. This also explains the difference between the expression of the applied field~\eqref{expression heff} in these articles. 

Our main result is the asymptotic stability of the skyrmion. Due to the lack of embedding $H^1 \hookrightarrow L^\infty$ in two space dimension, we chose to adapt our framework in order to consume a little extra regularity as possible. This prompted us to consider the following formalism. Let 
\begin{equation*}
    f:\R_+\times \R_+\rightarrow\R_+
\end{equation*} 
be a continuous non-zero function, and set for any $\eta\in H^2(\R^2),$\index{$N_f(\eta)$}
\begin{equation*}
    N(\eta)=f\big(\normHun{\eta},\norm{\eta}{H^2}\big).\index{$N_f(\eta)$}
\end{equation*}
Now we impose general conditions on $f$, under which the main result is true. Assume that there exists a universal constant $B_*>0$ such that all three following conditions are satisfied.
\begin{enumerate}[label=\roman*),ref={\roman*}]
    \item For any $(x,y),(x',y')\in \R_+\times \R_+$ such that $x\leq x'$ and $y\leq y'$, \label{1ere hypothese f}
    \[ f(x,y)\leq B_* f(x',y'), \]
    \item For any $x,y\in\R_ +$ and $\lambda>0, f(\lambda x,\lambda y)\leq \lambda B_*f(x,y)$,\label{2eme hypothese f}
    \item For any $\eta\in H^2(\R^2),\ \normLii{\eta}\leq B_*N_f(\eta)$. \label{3eme hypothese f}
\end{enumerate} 
(See below the theorem for comments on these conditions and meaningful examples of such $f$).

\begin{thm}\label{théorème: main stabilité asymptotique skyrmion}
Assume $h\geq h_0$ where $h_0$ is the constant given by Proposition~\ref{proposition: propriétés du skyrmion}. There exist positive constants $C_*,\tau_*$ and $\delta_*$ such that for any $m_0\in H_*^2(\R^2)$ satisfying 
\[ N_f(m_0-m_*)\leq \delta_*, \]
 then the solution $m$ associated with the initial data $m_0$ is globally defined, and there exists $x\in \q C^1(\R_+,\R^2)$ and $x_\ii\in\R^2$ such that\index{orthogonal decomposition}
\begin{equation*}
    m(t)=m_*\big( \cdot -x(t)\big) + \eta(t), 
\end{equation*}
with, for any $t\geq 0$, the bounds:
\begin{align*}
     \normHun{\eta(t)}+\normR{\dot{x}(t)}+ \normR{x(t)-x_\ii} & \le  C_* e^{-\tau_*t} \| m_0 - m_* \|_{H^1}, \\
\text{and} \qquad     \norm{\eta(t)}{H^2} & \le C_* e^{-\tau_*t}  \| m_0 - m_* \|_{H^2}.
\end{align*}

\end{thm}

The point of Theorem~\ref{théorème: main stabilité asymptotique skyrmion} is to provide asymptotic stability both at the $H^1$ (resp. at the $H^2$-level), but for perturbations which could large in $H^2$. Notice that asymptotic stability occurs with exponential rate, independent of the data.

We can hardly consider larger data spaces than $H^2_*$ since the local well-posedness itself holds in $H^2_*$.

Nevertheless, our result states that one has $H^2$ convergence (to a gauged skyrmion) even if the initial datum is near a skyrmion in merely $H^1$, with some control of the $H^2$ norm: the quantity $N_f$ quantifies this feature, which is reminiscent of parabolic equation.

We can give three significant examples of function $f$ that fulfill conditions~\eqref{1ere hypothese f}~\eqref{2eme hypothese f} and~\eqref{3eme hypothese f}, and for which Theorem~\ref{théorème: main stabilité asymptotique skyrmion} holds:
\begin{itemize}
    \item $f_1(x,y)=y$, which corresponds to the standard statement of asymptotic stability in $H_*^2$.
    \item $f_2(x,y)=\sqrt{xy}$ satisfies the condition~\eqref{3eme hypothese f} by the Gagliardo-Nirenberg inequality. 
    \item $f_3(x,y)=x\big(1+\sqrt{\log(1+\frac{y}{x})}\big)$, for which condition~\eqref{3eme hypothese f} holds true, due to the Brezis-Gallouët inequality.
\end{itemize}

With $f_2$, the initial perturbation need not be small in $H^2$, but inversely small in $H^2$ compared to the $H^1$-norm.
With $f_3$, the condition on the data is even less restrictive on the growth of the $H^2$-norm as the $H^1$-norm decreases. 

One could be tempted to recover the $H^2$ stability for small perturbations in the larger class $H^1$. This would correspond to a function $f_4(x,y)=x$, which does not satisfy condition~\eqref{3eme hypothese f} in dimension 2, and currently seems out of reach.

We lay emphasis on the fact that stability holds under general perturbation, and no further symmetry conditions are needed. This is consistent with the fact that the flow of~\eqref{llg} does not conserve the $O(2)$-equivariance (see~\cite{GustafWang}).

In other words, the solution associated with an initial data, taken sufficiently close to a skyrmion converges exponentially fast to a translated skyrmion. The method used in this article is reminiscent from~\cite{CotIgnat}, where the authors prove the asymptotic stability of domain walls under the flow of~\eqref{llg}, involving an extra time-dependent applied field. This article consists of several steps. Section~\ref{section: energy estimates} is devoted to the introduction of the main notation, some useful expansions, and eventually to estimates of the energy and its temporal derivatives, crucial for controlling the time-dependent perturbation. Then in Section~\ref{section: control of the perturbation}, we proceed to some spectral study of the linearized operator around the skyrmion and this leads us to the orthogonal decomposition of the solution. We emphasize that we require some $H^2$ control, but the variationnal properties of the dissipation lack $H^2$ coercivity (unlike what happens at $H^1$ regularity), see Remark \ref{rem:L_H^2}. Ultimately, we implement a bootstrap argument for smooth initial datum, and recover the asymptotic result by a regularizing argument.

\begin{notations}
Throughout the article, for any general functions $f,g$, we shall write\index{$\lesssim$} $f\lesssim g$ whenever there exists a constant $K>0$, only depending on the skyrmion $m_*$, such that $\left| f\right|\leq K \left|g\right|$.    
\end{notations}

\subsection{Acknowledgments}

The authors acknowledge partial support by the ANR project MOSICOF ANR-21-CE40-0004. R.C.
has benefitted from support provided by the University of Strasbourg Institute for Advanced Study
(USIAS) for a Fellowship within the French national program “Investment for the future” (IdEx-Unistra).

\section{Energy estimates}\label{section: energy estimates}

In this section, we proceed to some expansion of the energy, necessary to provide some control on perturbations of a skyrmion $m_*$.

\subsection{Expansion of the energy}

The energy $E$ can be expanded around $m_*$ as follows, for $\eta\in H^1(\R^2)$:
\begin{equation}\label{expression E(m_*+eta)=E(m_*)}
    E(m_*+\eta)=E(m_*)+\psLdeux{\delta E(m_*)}{\eta}+H(\eta,\eta),
\end{equation}
with\index{$\delta E$}
\begin{equation*}
    \delta E(m_*)=-\Delta m_* +2\grad\times m_* + h(m_*-e_3),
\end{equation*}

and\index{$H(\eta_1,\eta_2)$}
\begin{equation*}
    H(\eta_1,\eta_2)=\dfrac{1}{2}\int_{\R^2}\grad\eta_1.\grad\eta_2 + 2 \eta_1.\grad\times \eta_2 + h\eta_1.\eta_2,
\end{equation*}
satisfying $H(\eta_1,\eta_2)=H(\eta_2,\eta_1)$ and, for a constant $C_H$\index{$C_H$} only depending on $h$,
\begin{equation}\label{borne sur H(eta_1,eta_2)}
    H(\eta_1,\eta_2)\leq C_H \norm{\eta_1}{H^2}\normLdeux{\eta_2},
\end{equation}
for any $(\eta_1,\eta_2)\in H^2(\R^2)\times L^2(\R^2)$. In the sequel, we shall make crucial use of the identity
\begin{equation}\label{delta E(m_*+eta)= delta E(m_*) + Leta}
    \delta E(m_*+\eta)=\delta E(m_*) + L\eta,
\end{equation}

with the self-adjoint linearized operator\index{operator $L$} on $L^2(\R^2)$, given by
\begin{equation}\label{linéarisé L du gradient}
    L\eta=(-\Delta + 2\grad\!\!\times  + h )\eta ,
\end{equation}
with domain $H^2(\R^2)$. Recall from~\cite{LiMelcher} that the Euler-Lagrange equation associated with the geometric constraint $m\in\S^2$ reads
\begin{equation}
    \delta E(m_*) = \Lambda_*m_*\quad\text{with }\Lambda_*=|\grad m_*|^2 +2m_* . (\grad\times m_*) + h(1-\cos(\theta)).
\end{equation}

In the case of the skyrmion~\eqref{définition skyrmion}, the Lagrange multiplier $\Lambda_{*}$ reads
\begin{equation}\label{expression multiplicateur de Lagrange lambda m_*}
    \Lambda_{*}=(\theta')^2 + \dfrac{\sin^2(\theta) }{r^2} + 2\theta' + \dfrac{2\sin\theta\cos\theta}{r} + h(1-\cos\theta),
\end{equation}
and the property on the skyrmion number $Q(m_*)=-1$ is due to the boundary conditions satisfied by $\theta$.

\subsection{Moving frame}

We shall now proceed to a frame change and write the perturbation $\eta\in H^2(\R^2)$ along a new basis. This idea was introduced to deal with the stability of domain walls in ferromagnetic nanowires in~\cite{CarbLabb1} and then adapted in many related works such as~\cite{LiMelcher,CotIgnat}. In fact, we decompose $\eta$ onto a moving frame adapted to the skyrmion $m_*$. To do this, we set an orthogonal basis $(m_*,n_*,p_*)$ defined as
\begin{equation*}
    n_*=e_r\quad\text{and}\quad p_*=m_*\times n_* = \cos\theta\ e_\psi - \sin\theta\ e_3.
\end{equation*}

Then any $\eta$ can be decomposed as
\begin{equation}\label{décomposition de eta}
    \eta=\mu m_* + \nu n_* + \varpi p_*.
\end{equation}

A useful observation is that in view of the decomposition of $\eta$ on the orthogonal basis $(m_*,n_*,p_*)$, it is straightforward that for any integer $k$, if $\eta\in H^k(\R^2)$\index{$\lesssim$},
\begin{equation}\label{norm H^k eta et normH^k mu etc}
    \norm{\mu}{H^k}^2 +\norm{\nu}{H^k}^2 +\norm{\varpi}{H^k}^2\lesssim \norm{\eta}{H^k}^2\lesssim \norm{\mu}{H^k}^2 +\norm{\nu}{H^k}^2 +\norm{\varpi}{H^k}^2.
\end{equation}

In terms of the moving basis, the polar coordinates read
\begin{align*}
    &e_r=n_*\\
    & e_\psi=\sin\theta \ m_* + \cos\theta \ p_*\\
    &e_3=\cos\theta\  m_*-\sin\theta\  p_*.
\end{align*}

Imposing the geometrical constraint $m_*+\eta\in\mathbb{S}^2$ has important consequences, such as the fact that 
\begin{equation}\label{mu en fonction de eta^2}
    \mu = \eta.m_* = -\frac{|\eta|^2}{2}=-\dfrac{\mu^2+\nu^2+\varpi^2}{2}.
\end{equation}

Furthermore, one can compute
    \begin{align*}
    \partial_r m_* = \theta' p_*,\quad &\partial_\psi m_* = -\sin\theta\ n_*,\quad \\
     \partial_r n_*=0,\quad & \partial_\psi n_*=\sin\theta\ m_* + \cos\theta\ p_*,\\
    \partial_r p_* = -\theta' m_*, \quad &\partial_\psi p_* = -\cos\theta\ n_*,
    \end{align*}
    as well as
    \begin{align*}
        \partial_1 m_* = \dfrac{\sin\psi\sin\theta}{r}n_*+\cos\psi\theta'p_*,\quad &\partial_2 m _*=-\dfrac{\cos\psi\sin\theta}{r}n_*+\sin\psi\theta'p_*,\\
        \partial_1 n_* = -\dfrac{\sin\psi}{r}\big(\sin\theta \ m_* + \cos\theta p_*),\quad &\partial_2 n _0=\dfrac{\cos\psi}{r}\big(\sin\theta m_* + \cos\theta\ p_*),\\
        \partial_1 p_* = -\cos\psi \theta'm_*+\dfrac{\sin\psi}{r}\cos\theta \ n_*,\quad &\partial_2 p_*=-\sin\psi\theta' m_*-\dfrac{\cos\psi}{r}\cos\theta \ n_*.\\
    \end{align*}

We now write the general form of the operator $L$ on the orthonormal basis $(m_*,n_*,p_*)$.

\begin{lem}\label{lemme: expression Leta}
For any $\eta\in H^2(\R^2)$ decomposing as $\eta = \mu m_* + \nu n_* +\varpi p_*$, then $\mu,\nu,\varpi\in H^2(\R^2)$ and
    \begin{align*}
        L\eta=&(-\Delta\mu+M_\varpi\varpi + M_{\varpi,r}\partial_r\varpi+M_\mu\mu + M_{\nu,\psi}\partial_\psi\nu)m_* \\
        &+ (-\Delta \nu + N_\nu\nu - M_{\nu,\psi}\partial_\psi \mu -P_{\nu,\psi}\partial_\psi\varpi)n_* \\
    & + (-\Delta \varpi + P_\varpi\varpi + P_\mu \mu -M_{\varpi,r}\partial_r\mu + P_{\nu,\psi}\partial_\psi\nu)p_*, 
    \end{align*}
    with
    \begin{align*}
    &M_\varpi=\dfrac{2\cos\theta\sin\theta}{r^2}+\dfrac{2\cos(2\theta)}{r}+h\sin\theta ,\\
    &M_\mu=\theta'^2+\frac{\sin^2\theta}{r^2}+2\theta'+\dfrac{2\cos\theta\sin\theta}{r}+h,\\
    &M_{\varpi,r} = 2\theta'+2,\\
    & M_{\nu,\psi}=-\dfrac{2}{r}\left(\dfrac{\sin\theta}{r}+\cos\theta\right),\\
    & N_\nu = \dfrac{1}{r^2}+h ,\\
    & P_\varpi = 2\theta'-\dfrac{2\sin\theta\cos\theta}{r}+h+(\theta')^2+\dfrac{\cos^2(\theta)}{r^2},\\
    & P_\mu = -h\sin\theta,\\
    & P_{\nu,\psi}=\dfrac{2}{r}\left(\sin\theta-\dfrac{\cos\theta}{r}\right).
\end{align*}
\index{$P_\varpi$}
\index{$N_\nu$}
\index{$P_{\nu,\psi}$}
\end{lem}

\begin{proof}[Proof of Lemma~\ref{lemme: expression Leta}]

By~\eqref{norm H^k eta et normH^k mu etc}, for $\eta\in H^2(\R^2)$, we are able to give a proper sense to all following quantities. We compute
    \begin{align*}
    &\partial_r\eta =(\partial_r\mu-\theta'\varpi)m_*  + \partial_r\nu\ n_* + (\mu\theta'+\partial_r\varpi)p_*,
\end{align*}
and
\begin{align*}
    \partial_\psi\eta=(\partial_\psi \mu +\nu\sin\theta)m_* + (\partial_\psi\nu -\mu\sin\theta-\varpi\cos\theta)n_* + (\nu\cos\theta+\partial_\psi\varpi)p_*,
\end{align*}

as well as
\begin{align*}
    \partial_r^2\eta = \big[ \partial_r^2\mu - \theta'' \varpi-&(2\partial_r\varpi+\mu\theta')\theta'\big]m_*  + \partial_r^2\nu\  n_* \\
    &+\big[(2\partial_r\mu - \theta'\varpi)\theta' +\mu\theta'' + \partial_r^2\varpi\big]p_*
    ,
\end{align*}
and
\begin{align*}
    \partial_\psi^2\eta =& \big[ \partial_\psi^2\mu + (2\partial_\psi\nu-\mu\sin\theta-\varpi\cos\theta)\sin\theta\big]m_* \\
    &+\big[\partial_\psi^2\nu -2\partial_\psi\mu \sin\theta - \nu - 2 \partial_\psi\varpi\cos\theta\big]n_* \\
    &+\big[\partial_\psi^2\varpi + (2\partial_\psi\nu - \mu\sin\theta-\varpi\cos\theta)\cos\theta\big]p_*.
\end{align*}

This leads to 
\begin{align*}
    \Delta\eta  = &\bigg[\Delta\mu
-\Big(\theta'^2+\frac{\sin^2\theta}{r^2}\Big)\mu
-\Big(\theta''+\frac{1}{r}\theta'+\frac{1}{r^2}\sin\theta\cos\theta\Big)\varpi\\
&+ \frac{2}{r^2}\sin\theta\,\partial_\psi\nu
-2\theta'\,\partial_r\varpi\bigg]m_*\\
&+ \bigg[\Delta\nu
-\frac{1}{r^2}\nu
-\frac{2}{r^2}\sin\theta\,\partial_\psi\mu
-\frac{2}{r^2}\cos\theta\,\partial_\psi\varpi\bigg]n_*\\
& + \bigg[\Delta\varpi
-\Big(\theta'^2+\frac{\cos^2\theta}{r^2}\Big)\varpi
+\Big(\theta''+\frac{1}{r}\theta'-\frac{1}{r^2}\sin\theta\cos\theta\Big)\mu \\
&+2\theta'\,\partial_r\mu
+\frac{2}{r^2}\cos\theta\,\partial_\psi\nu\bigg]p_*.
\end{align*}

On the other hand, we compute
\begin{align*}
    \grad\times \eta  =& \bigg[\partial_r\varpi + \mu\,\theta'
   + \dfrac{\cos\theta}{r}\big(\mu\sin\theta+\varpi\cos\theta-\partial_\psi\nu\big)\bigg]m_* \\
   &+ \bigg[\dfrac{1}{r}\big(\partial_\psi\mu\cos\theta-\partial_\psi\varpi\sin\theta\big)\bigg]n_*\\
   & +\bigg[-\partial_r\mu + \varpi\,\theta'
   - \dfrac{\sin\theta}{r}\big(\mu\sin\theta+\varpi\cos\theta-\partial_\psi\nu\big)\bigg]p_*,
\end{align*}

and
\begin{equation*}
    h(\eta-e_3)= h\big[(\mu-\cos\theta)m_* +\nu n_* +(\varpi+\sin\theta)p_*\big].
\end{equation*}

To conclude, we write
\begin{align*}
    L\eta=&-\Delta\eta +2\grad\times \eta + h\eta\\
    =&\bigg[-\Delta\mu + \left((\theta')^2 + \dfrac{\sin^2(\theta)}{r^2}+2\theta' + \dfrac{2\cos\theta\sin\theta}{r}+h\right)\mu\\
    &+\left(\theta'' + \dfrac{\theta'}{r}+\dfrac{\cos\theta\sin\theta}{r^2}+2\dfrac{\cos^2(\theta)}{r}\right)\varpi \\
    & +(2\theta'+2)\partial_r\varpi-\dfrac{2}{r}\left(\dfrac{\sin\theta}{r}+\cos\theta\right)\partial_\psi\nu\bigg]m_*\\
    &+\bigg[-\Delta\nu+\Big(\dfrac{1}{r^2}+h\Big)\nu+\dfrac{2}{r}\Big(\dfrac{\sin\theta}{r}+\cos\theta\Big)\partial_\psi\mu \\
    &+ \dfrac{2}{r}\Big(\dfrac{\cos\theta}{r}-\sin\theta\Big)\partial_\psi\varpi\bigg]n_*\\
    &+\bigg[-\Delta\varpi + \Big(2\theta'-\dfrac{2\sin\theta\cos\theta}{r}+h+(\theta')^2+\dfrac{\cos^2(\theta)}{r^2}\Big)\varpi\\
    &-\left(\theta'' + \dfrac{\theta'}{r}-\dfrac{\cos\theta\sin\theta}{r^2}+\dfrac{2\sin^2(\theta)}{r}\right)\mu\\
    &-(2\theta'+2)\partial_r\mu+\dfrac{2}{r}(\sin\theta-\dfrac{\cos\theta}{r})\partial_\psi\nu\bigg]p_*\\
    =&\ (-\Delta\mu+M_\varpi\varpi + M_{\varpi,r}\partial_r\varpi+M_\mu\mu + M_{\nu,\psi}\partial_\psi\nu)m_*  \\
    & + (-\Delta \nu + N_\nu\nu - M_{\nu,\psi}\partial_\psi \mu -P_{\nu,\psi}\partial_\psi\varpi)n_*\\
    & + (-\Delta \varpi + P_\varpi\varpi +P_\mu \mu -M_{\varpi,r}\partial_r\mu + P_{\nu,\psi}\partial_\psi\nu)p_* . \qedhere
\end{align*}
\end{proof}

As in~\cite{CotIgnat}, we now introduce useful notations to simplify the study.

\begin{notations}
For integers $k\geq 0$ and $l\geq 1$ and a vector valued function $\eta$. We shall write that $u= \mathcal{O}_k^l(\eta)$ if the function $u$ is an homogeneous polynomial of total degree $l$ in the components of $\eta$ and their derivatives up to order $k$ such that the total number of derivatives in each term is at most $k$, and the coefficients are $W^{2,\ii}(\R^2)$ functions (this notation is reminiscent from~\cite{CotIgnat} and we refer to Definition 4.5 in the latter article for further examples).
\end{notations}

Throughout the article, we refer to Appendix~\ref{appendix: estimates for high order terms} for useful estimates regarding terms expressed in terms of $\int_{\R^2}\mathcal{O}_k^l(\eta)$\index{$\mathcal{O}_k^l(\eta)$}. In particular, using~\eqref{mu en fonction de eta^2}, we introduce a crucial example, which we often use throughout the article, that is
\begin{equation}\label{mu = quadratique en nu}
    \mu =-\dfrac{\nu^2+\varpi^2}{2} + \mathcal{O}_0^4(\eta).
\end{equation}

\subsection{A first estimate on the energy}

\begin{lem}\label{lemme: equivalence norm H^1 eta et (pi,nu)}
    There exists two positive constants $C_0$ and $\delta_0$ such that for $\eta\in H^2(\R^2)$ decomposing according to~\eqref{décomposition de eta} and such that $N_f(\eta)\leq \delta_0$,\index{$N_f(\eta)$} we have
    \begin{equation*}
        \dfrac{1}{C_0}\norm{(\nu,\varpi)}{H^2}^2\leq \norm{\eta}{H^2}^2\leq C_0\norm{(\nu,\varpi)}{H^2}^2,
    \end{equation*}

    and
    \begin{equation*}
        \dfrac{1}{C_0}\norm{(\nu,\varpi)}{H^1}^2\leq \norm{\eta}{H^1}^2\leq C_0\norm{(\nu,\varpi)}{H^1}^2.
    \end{equation*}
\end{lem}

\begin{proof}[Proof of Lemma~\ref{lemme: equivalence norm H^1 eta et (pi,nu)}]
By~\eqref{norm H^k eta et normH^k mu etc}, both desired left-hand side inequalities are valid.
Moreover,\index{$\lesssim$} by~\eqref{mu en fonction de eta^2}, and the two-dimensional Sobolev embeddings and the Gagliardo-Nirenberg inequality,
\begin{align*}
    \norm{\mu}{H^2}&\lesssim \normLdeux{|\eta|^2} + \normLdeux{|\grad^2\eta||\eta|}+\normLdeux{|\grad\eta|^2}\\
    &\lesssim \norm{\eta}{L^4}^2 + \normLdeux{\grad^2\eta}\normLii{\eta} + \norm{\grad\eta}{L^4}^2\\
    &\lesssim \norm{\eta}{H^2}\Big(\normHun{\eta} + \normLii{\eta}\Big).
\end{align*}

In a same manner, we have
\begin{align*}
    \norm{\mu}{H^1}&\lesssim \normLdeux{|\eta|^2} + \normLdeux{|\grad\eta||\eta|}\\
    &\lesssim \norm{\eta}{H^1}^2 + \norm{\grad\eta}{L^2}\norm{\eta}{L^\ii}\lesssim  \norm{\eta}{H^1}\big(\norm{\eta}{H^1} + \normLii{\eta}\big).
\end{align*}

By Proposition~\ref{proposition: propriétés de f}, we have $\normHun{\eta} + \normLii{\eta}\leq CN_f(\eta) \le C \delta_0$, and we obtain the first desired estimate in Lemma~\ref{lemme: equivalence norm H^1 eta et (pi,nu)}.
The second desired estimate follows by taking $\delta_0$ small enough.
\end{proof}

\begin{lem}\label{lemme: controle estimation E(m_*+eta)-E(m_*)}
There exist two positive constants $C_1$ and $\delta_1$ such that for any $\eta\in H^2(\R^2)$ that decomposes according to~\eqref{décomposition de eta} and satisfying $N_f(\eta)\leq \delta_1$, we have
\begin{equation*}
    E(m_*+\eta) -E(m_*)\leq C_ 1\normHun{\eta}^2 
\end{equation*}
and
\begin{equation*}
     - C_1N_f(\eta)\norm{\eta}{H^1}^{2}\leq E(m_*+\eta) -E(m_*)-\dfrac{1}{2}\psLdeuxLdeux{\q L(\nu,\varpi)}{(\nu,\varpi)}
\end{equation*}
where the operator $\q L$ has a domain $\dom(\q L)=H^2(\R^2)^2$ and is given by \index{operator $\q L$}\index{operator $\q L_\nu$}\index{operator $\q L_\varpi$}
\begin{gather} \label{def:qL}
    \q L=\begin{pmatrix}
             \q L_\nu & S\\
             S^* & \q L_\varpi
         \end{pmatrix}\text{ where}\quad S=-P_{\nu,\psi}\partial_\psi, \\
         \q L_\nu=-\Delta + N_\nu-\Lambda_{*} \quad \text{and} \quad
         \q L_\varpi=-\Delta + P_\varpi - \Lambda_{*}. \label{def:qLnu_pi}
\end{gather}
\index{$P_\varpi$}
\index{$N_\nu$}
\index{$P_{\nu,\psi}$}
\end{lem}

\begin{rem}\label{remarque: différence entre les opérateurs linéarizés L}
With the above notation, we observe that $L$ and $\q L$ are related by the formula 
    \begin{equation*}
        \big(L-\Lambda_{*}\big)\eta = \q L_\mu(\nu,\varpi) m_* + (\q L_\nu\nu + S\varpi)n_* + (\q L_\varpi\varpi + S^*\nu)p_* + \mathcal{O}_1^2(\eta),
    \end{equation*}
with $\q L_\mu(\nu,\varpi)=-\Delta\mu+M_\varpi\varpi + M_{\varpi,r}\partial_r\varpi+\big(M_\mu-\Lambda_{*}\big)\mu + M_{\nu,\psi}\partial_\psi\nu$. The operator $\q L$ corresponds to the linearization of the energy around a skyrmion in both $\nu$ and $\varpi$ variables. Due to the geometry of the sphere, we recall that we have~\eqref{mu = quadratique en nu} and in view of~\eqref{delta E(m_*+eta)= delta E(m_*) + Leta}, $L-\Lambda_{*}$ in fact corresponds to the linearized operator around a skyrmion in the original setting.
\end{rem}

We now prove Lemma~\ref{lemme: controle estimation E(m_*+eta)-E(m_*)}.

\begin{proof}[Proof of Lemma~\ref{lemme: controle estimation E(m_*+eta)-E(m_*)}]
Using~\eqref{expression E(m_*+eta)=E(m_*)} and then Lemma~\ref{lemme: expression Leta}, we derive
\begin{align*}
    &E(m_*  +\eta)-E(m_*)  =\int_{\R^2}\big( \delta E(m_*).\eta\big) + H(\eta,\eta)\\
    &=\dfrac{1}{2}\int_{\R^2} |\grad\eta|^2 + 2\eta.\grad\times \eta + (h-\Lambda_{*})|\eta|^2\\
    &=\dfrac{1}{2}\Big(\psLdeux{\q L_\nu(\nu)}{\nu}+\psLdeux{\q L_\varpi(\varpi)}{\varpi} \Big)+ 2 \int_{\R^2}P_{\nu,\psi}\varpi \partial_\psi\nu \\
    &+\int_{\R^2} \Big(\dfrac{|\grad \mu|^2 }{2}- \dfrac{(M_\mu-\Lambda_*)\mu^2}{2}+\dfrac{(M_\varpi+P_\mu) \varpi \mu}{2}+ M_{\varpi,r}\partial_r\varpi\mu + M_{\nu,\psi}\partial_\psi\nu\mu\Big).
\end{align*}

First, $|\grad\mu|^2 = \mathcal{O}_2^4(\eta)$. Then, the asymptotics given in Proposition~\ref{proposition: propriétés du skyrmion} and the expression of the Lagrange multiplier in~\eqref{expression multiplicateur de Lagrange lambda m_*} imply that $\Lambda_*,M_{\varpi,r},M_\mu $ and $P_\mu$ as functions of the variable $x\in\R^2$ (or $(r,\psi)\in [0,+\ii)\times [0,2\pi)$). Therefore, all terms attached to such functions in the previous integrand can be viewed as $\mathcal{O}_1^4(\eta)$ or $\mathcal{O}_1^3(\eta)$. Moreover, referring to the cylindrical formulation of the gradient on $\R^2$, we obtain $M_{\nu,\psi}\partial_\psi\nu=\mathcal{O}_1^2(\eta)$. At last, we use that $x\mapsto |x|M_\varpi$ is bounded on $\R^2$, and we write 
\begin{equation*}
    \left|\int_{\R^2}M_\varpi\varpi\mu\right|\leq \normLii{|x|M_\varpi}\normLdeux{\dfrac{\varpi}{|x|}}\normLdeux{\mu}\lesssim \normLdeux{\dfrac{\varpi}{|x|}}\normHun{\eta}^2.
\end{equation*}

Separating the case $|x|\geq 1$ and $|x|\leq 1$ in $\normLdeux{\frac{\varpi}{|x|}}$, then using Holder's inequality and Gagliardo-Nirenberg inequality, we obtain $\normLdeux{\frac{\varpi}{|x|}}\lesssim \normHun{\varpi}$.

Combining all previous estimates, and by Lemma~\ref{lemme: controle int O^l_1 et O^l_2}, we have 
\begin{equation*}
    \left|\int_{\R^2}M_\varpi\varpi\mu\right|+\Big|\int_{\R^2}\mathcal{O}_1^3(\eta)\Big|+\Big|\int_{\R^2}\mathcal{O}_2^4(\eta)\Big|\lesssim \normHun{\eta}^2\Big(\normHun{\eta} + \normLii{\eta}^2\Big).
\end{equation*}
By Proposition~\ref{proposition: propriétés de f}, $\normHun{\eta} + \normLii{\eta}\leq CN_f(\eta)$, which leads to the desired estimate.\end{proof}

\subsection{Dissipation term}\label{section: dissipation term}

Assume that the solution $m(t)$ is well-defined and smooth enough so that $E\big(m(t)\big)$ is differentiable with respect to the variable $t$, as stated in Theorem~\ref{théorème: LWP of the Cauchy problem}. These statements will be checked later and we refer to Appendix~\ref{appendix: LWP} for more details. We anticipate the former facts and formally differentiate $E\big(m(t)\big)$ and obtain
\begin{align}
    \dfrac{d}{dt}\big(E(m(t)\big)&=\int_{\R^2}\delta E\big(m(t)\big).\partial_tm(t)\label{dérivée de E(m(t)) non rigoureuse}\\
    &=-\int_{\R^2}\heff(m).\Big( m\times \heff(m) - \alpha m\times\big(m\times \heff(m)\big)\Big)\notag\\
    &=-\alpha\int_{\R^2} \big(|\delta E(m)|^2-(\delta E(m).m)^2\big)\notag.
\end{align}

This now motivates the rigorous definition. For $m\in H_*^2$, $\delta E(m)\in L^2(\R^2)$ and $m\in L^\ii(\R^2)$, we then define\index{$Diss_1$}
\begin{equation}
    Diss_1(m)=-\alpha\int_{\R^2} \big(|\delta E(m)|^2-(\delta E(m).m)^2\big),
\end{equation}

and we have the following result.
\begin{lem}\label{lemme: controle perturbation niveau H^1}
There exist two positive constants $C_2$ and $\delta_2$, such that for any $\eta\in H^2(\R^2)$ satisfying $N_f(\eta)\leq \delta_2$, \index{$N_f(\eta)$}we have
\begin{equation*}
    Diss_1(m) \leq -\alpha \normLdeuxLdeux{\q L(\nu,\varpi)}^2 + C_2\delta_2\norm{\eta}{H^1}^2 .
\end{equation*}
\end{lem}

\begin{proof}[Proof of Lemma~\ref{lemme: controle perturbation niveau H^1}]
Since $\eta\in H^2(\R^2)$, in view of the properties of the skyrmion $m_*$, we deduce that $m=m_* + \eta\in H^2_*$. Then, using the expression of $L\eta$ in the proof of Lemma~\ref{lemme: expression Leta} and the notations introduced in Lemma~\ref{lemme: controle estimation E(m_*+eta)-E(m_*)}, we compute   
\begin{align*}
     Diss_1(m) &=-\alpha\int |\delta E(m_*+\eta)|^2 - \big((m_*+\eta).\delta E(m_*+\eta)\big)^2\\
    &=-\alpha\int_{\R^2} \bigg(|L\eta|^2 - \big(m_*.L\eta\big)^2 - 2\big(\Lambda_{*}\big)^2\eta.m_* - 2\Lambda_{*}\eta. L\eta\bigg) \\
    &\quad \quad \quad +\int_{\R^2}\mathcal{O}_2^3(\eta) +\int_{\R^2}\mathcal{O}_2^4(\eta)\\
    &=-\alpha\normLdeux{\q L_\nu \nu-P_{\nu,\psi}\partial_\psi\varpi}^2 -\alpha\normLdeux{\q L_\varpi \varpi+P_{\nu,\psi}\partial_\psi\nu}^2 \\
    &\quad \quad \quad +\int_{\R^2}\mathcal{O}_2^3(\eta) +\int_{\R^2}\mathcal{O}_2^4(\eta)\\
    &=-\alpha\normLdeuxLdeux{\q L(\nu,\varpi)}^2 +\int_{\R^2} \mathcal{O}_2^3(\eta) +\int_{\R^2}\mathcal{O}_2^4(\eta).
\end{align*}

By Lemma~\ref{lemme: controle int O^l_1 et O^l_2}, we obtain 
\begin{align*}
    \Big|Diss_1(m) + \alpha \normLdeuxLdeux{\q L(\nu,\varpi)}^2\Big| \leq \normHun{\eta}^2 \Big(\normLii{\eta} + \normLii{\eta}^2\Big).
\end{align*}

Using~\eqref{3eme hypothese f}, and taking $\delta_2$ small enough, we deduce that 
\begin{equation*}
    \Big|Diss_1(m) + \alpha \normLdeuxLdeux{\q L(\nu,\varpi)}^2\Big| \leq  C_2\delta_2\norm{\eta}{H^1}^2. \qedhere
\end{equation*}

\end{proof}

\subsection{Higher order term}\label{section: higher order term}

Contrary to the one-dimensional case, as dealt with in~\cite{CotIgnat}, the control provided in Lemma~\ref{lemme: controle perturbation niveau H^1} will not be sufficient to prove the stability result. This section is devoted a new quantity, that eventually provides good estimates for the second order derivatives. To overcome this structural problem, we introduce an higher order dissipative term. Once more, we anticipate\footnote{See Subsection~\ref{section: Bootstrap argument for a smooth initial data} for more details.} that the solution $m(t)$ is smooth enough and decomposes as $m(t)=m_*\big(.-x(t)\big) + \eta(t)$ with a differentiable function $x$, so that 
\[ t\mapsto \normLdeux{L\eta(t)}^2=\normLdeux{\delta E\big(m(t)\big)-\delta E\big(m_*(.-x(t))\big)}^2 \]
is differentiable (recall the definition of $L$ \eqref{delta E(m_*+eta)= delta E(m_*) + Leta}-\eqref{linéarisé L du gradient}).
At least formally, and using~\eqref{borne sur H(eta_1,eta_2)}, we derive\index{$H(\eta_1,\eta_2)$}
\begin{align}
    \dfrac{1}{2}\dfrac{d}{dt} \| L\eta \|_{L^2}^2 & =H\Big(\partial_t m(t)+\dot{x}(t)\partial_x m_*\big(.-x(t)\big),L\eta(t) \Big)\notag\\
    &\leq\int_{\R^2} \Big(m(t)\times \heff(m(t))-\alpha m\times \big(m(t)\times\heff(m(t))\big)\Big).L^2\eta(t)\label{dérivée ordre deux non rigoureuse} \\
    &\quad\quad +C_H|\dot{x}(t)|\norm{m_*}{Z^3}\normLdeux{L\eta(t)}.\notag
\end{align}

For simplicity, we can cast the dependence upon the translation $x(t)$ for the next observation. For a magnetization that decomposes as $m=m_*+\eta$ with $\eta\in H^3(\R^2)$, we can give a proper sense to
\begin{equation*}
    Diss_2(m):= \int_{\R^2} \Big(m\times \heff(m)-\alpha m\times \big(m\times\heff(m)\big)\Big).L^2\eta.
\end{equation*}

Furthermore, we can also deduce a control of $Diss_2(m)$ in terms of the Sobolev norms up to the third order. These claims are summarized in the following lemma.

\begin{lem}\label{lemme: controle au niveau H^2 perturbation}
    There exist two positive constants $C_3>0$ and $\delta_3>0$ such that for any $\eta\in H^3(\R^2)$, $Diss_2(m_* + \eta)$ is well-defined. In addition, if $N_f(\eta)\leq \delta_3$, we have 
    \begin{equation*}
        Diss_2(m)  + \dfrac{\alpha}{2}\norm{\eta}{\dot{H}^3}^2\leq C_3\big(\norm{\eta}{H^2}^2 +\norm{\eta}{H^2}^4\big) .
    \end{equation*}
\end{lem}

\begin{proof}[Proof of Lemma~\ref{lemme: controle au niveau H^2 perturbation}]
    
To deal with general $\eta\in H^3(\R^2)$, we invoke a density argument. Take a smooth sequence of functions $(\eta^n)$ with compact support tending to $\eta$ in $H^3$, then we establish the desired estimates for these smooth perturbations. Since the next quantities eventually do not involve more than third order derivatives, we shall recover the desired estimates for $\eta\in H^3(\R^2)$ by letting $n$ tend to $+\ii$. First, take $\eta\in C^\ii_0(\R^2)$ and split the main quantity into two parts,
\begin{align*}
    Diss_2(m)=\mathcal{I} + \mathcal{J},
\end{align*}
with
\begin{align*}
    \mathcal{I}&= \int_{\R^2} m\times\heff(m) .\Big((-\Delta + 2 \grad\times  + h )^2\eta\Big),
\end{align*}
and
\begin{align*}
    \mathcal{J}& = -\alpha\int_{\R^2} m\times (m\times \heff(m)).\Big((-\Delta + 2 \grad\times  + h )^2\eta\Big).
\end{align*}

Using that $m_*$ is a critical point, we shall linearize and obtain
\begin{align}\label{expression de I}
    \mathcal{I} = \int_{\R^2} \Big(\eta \times \heff(m_*)-m_*\times L\eta \Big).L^2\eta -\int_{\R^2} \eta\times L\eta .L^2\eta, 
\end{align}

and\begin{align}
    \mathcal{J}=& - \alpha\int_{\R^2} \Big(m_*\times ( \eta\times \heff(m_*))-m_*\times (m_*\times L\eta)  \Big).L^2\eta\label{expression de J} \\
    &-\alpha\int_{\R^2} \Big(\eta\times (\eta\times \heff(m_*)) - m_*\times ( \eta\times L\eta)-\eta\times (m_*\times L\eta)\Big). L^2\eta \notag\\
    &+\alpha \int_{\R^2} \eta\times (\eta\times L\eta).L^2\eta\notag
\end{align}

We set $L_0=-\Delta$\index{$L_0$} and we claim the following lemma.
\begin{lem}\label{lemme: int(atimes b).Lc = int (Latimes b).c +int (atimes L_0 b).c + G(a,b,c)}
For any $a,b,c\in H^2(\R^2)$, we have the identity
    \begin{equation*}
        \int_{\R^2}(a\times b).Lc = \int_{\R^2} (La\times b).c +\int_{\R^2} (a\times L_0 b).c + G(a,b,c),
    \end{equation*}
    where
    \begin{equation*}
         G(a,b,c)=-2\sum_{j=1}^3\int_{\R^2} (\partial_j a\times \partial_jb).c + 2\int_{\R^2}\Big(\grad\times (a\times b)-\big((\grad\times a) \times b\big).c
    \end{equation*}
\end{lem}

In view of Lemma~\ref{lemme: int(atimes b).Lc = int (Latimes b).c +int (atimes L_0 b).c + G(a,b,c)}, we deduce the following claims.

\begin{claim}\label{claim: controle de I}
For any $\varepsilon>0$, there exist a constant $K_1>0$ independant of $\varepsilon$ and $K_\varepsilon>0$ (depending on $\varepsilon$) such that
    \begin{equation}
    \big|\mathcal{I}\big|\leq K_1\norm{\eta}{\dot{H}^3}^2\Big(\varepsilon + \normLii{\eta} + \normHun{\eta}^\frac{2}{3}\Big) + K_\varepsilon\Big(\norm{\eta}{H^2}^2 + \norm{\eta}{H^2}^4\Big).
\end{equation}
\end{claim}

\begin{proof}[Proof of Claim~\ref{claim: controle de I}]
We have $\big(\eta\times\heff(m_*)\big).L^2\eta =\mathcal{O}_0^1(\eta)\mathcal{O}_4^1(\eta)= \mathcal{O}_4^2(\eta)$, hence 
\begin{equation*}
    \int_{\R^2}\big(\eta\times \heff(m_*)\big).L^2\eta = \int_{\R^2}\mathcal{O}_4^2(\eta) . 
\end{equation*}

On the other hand, using Lemma~\ref{lemme: int(atimes b).Lc = int (Latimes b).c +int (atimes L_0 b).c + G(a,b,c)}, we deduce 
\begin{align*}
    -\int_{\R^2} \big(m_*\times L\eta \big).L^2\eta =& -\dfrac{1}{2}\int_{\R^2} \big(m_*\times L\eta \big).L^2\eta - \dfrac{1}{2}\int_{\R^2} \big(Lm_*\times L\eta\big).L\eta  \\
    & - \dfrac{1}{2}\left( \int_{\R^2} \big(m_*\times L_0L\eta \big).L\eta +G(m_*,L\eta,L\eta)\right).
\end{align*}

Due to the fact that 
\begin{align*}
    -\dfrac{1}{2}\int_{\R^2} \big(m_*\times L\eta \big).L^2\eta =&  -\dfrac{1}{2}\int_{\R^2}\big(m_*\times L_0\eta \big).(2L_0N+N^2)\eta\\
    &-\dfrac{1}{2}\int_{\R^2} \big(m_*\times L_0\eta \big).L_0^2\eta -\dfrac{1}{2}\int_{\R^2} \big(m_*\times N\eta \big).L^2\eta
\end{align*}
and that
\begin{align*}
    -\dfrac{1}{2}\int_{\R^2} \big(m_*\times L_0L\eta \big).L\eta =&-\dfrac{1}{2} \int_{\R^2} \Big(\big(m_*\times L_0^2\eta \big).L_0\eta + \big(m_*\times L_0N\eta \big).L\eta \Big)  \\
    &-\dfrac{1}{2} \int_{\R^2}\big(m_*\times L_0^2\eta \big).N\eta\\
    = &\dfrac{1}{2}\int_{\R^2} \big(m_*\times L_0\eta \big).L_0^2\eta +\int_{\R^2} \mathcal{O}_5^2(\eta),
\end{align*}

we infer
\begin{equation*}
    -\int_{\R^2} \big(m_*\times L\eta \big).L^2\eta = G(m_*,L\eta,L\eta) + \int_{\R^2} \mathcal{O}_5^2(\eta) = \int_{\R^2} \mathcal{O}_5^2(\eta). 
\end{equation*}

It remains to deal with the trilinear term in~\eqref{expression de I}. We deduce\index{$L_0$}
\begin{align*}
    -\int_{\R^2} \big(\eta\times L\eta\big).L^2\eta & = \int_{\R^2} \mathcal{O}_6^3(\eta).
\end{align*}

Let $\varepsilon>0$. relying on Lemma~\ref{lemme: controle int O^l_3 et O^l_4}, we use classical bounds of the type $|ab|\leq \frac{\varepsilon}{p} a^p + \frac{1}{q\varepsilon} b^q$ with $1=\frac{1}{p} + \frac{1}{q}$, so that we deduce the existence of a constant $C_\varepsilon$ such that
\begin{equation*}
    \left|\int_{\R^2} \mathcal{O}_4^2(\eta)\right|+\left|\int_{\R^2} \mathcal{O}_5^2(\eta)\right|\lesssim \varepsilon\norm{\eta}{\dot{H}^3}^2 + C_\varepsilon\norm{\eta}{H^2}^2,
\end{equation*}
and in addition, using the above Hölder inequality with $p= \frac{4}{3}$ implies
\begin{equation*}
    \left|\int_{\R^2} \mathcal{O}_6^3(\eta)\right|\lesssim \norm{\eta}{\dot{H}^3}^2\Big(\varepsilon + \normLii{\eta}+\normHun{\eta}^\frac{2}{3}\Big) + C_\varepsilon\norm{\eta}{H^2}^4,
\end{equation*}
which concludes the proof of this claim.
\end{proof}

\begin{claim}\label{claim: controle de J + alpha norm H^3 eta}
For any $\varepsilon>0$, there exist a constant $\widetilde{K}_1>0$ independant of $\varepsilon$ and $\widetilde{K}_\varepsilon>0$ (depending on $\varepsilon$) such that
    \begin{align*}
    \mathcal{J}+\alpha\norm{\eta}{\dot{H}^3}^2\leq \widetilde{K}_1\norm{\eta}{\dot{H}^3}^2\Big(\varepsilon + &\normLii{\eta}^2 + \normHun{\eta}^2 + \normLii{\eta}\normHun{\eta}^\frac{2}{3}\Big) \\
    &+ \widetilde{K}_\varepsilon\Big(\norm{\eta}{H^2}^2 + \norm{\eta}{H^2}^4\Big).
\end{align*}
\end{claim}

\begin{proof}[Proof of Claim~\ref{claim: controle de J + alpha norm H^3 eta}]
Starting with the terms bilinear in $\eta$, we deal with the analysis of the higher order term that reads, using that $N$ and $L_0$ commute (see it in Fourier variables),
\begin{align*}
    \alpha \int_{\R^2} m_*\times (m_*\times L\eta).L^2\eta = \mathcal{J}_1  + \mathcal{J}_2,
\end{align*}
with
\begin{align*}
    & \mathcal{J}_1 = \alpha\int_{\R^2}   m_*\times (m_*\times L_0\eta).L_0^2\eta,\\
    &\mathcal{J}_2 = \alpha\int_{\R^2}  m_*\times (m_*\times L_0\eta).(2L_0N+N^2)\eta + \alpha \int_{\R^2}  m_*\times (m_*\times N\eta).L^2\eta.
\end{align*}

On the one hand, we have 
\begin{align*}
    \mathcal{J}_1 =  &\ \alpha\int_{\R^2} (m_*.L_0\eta )(m_*.L_0^2\eta)-\alpha\int_{\R^2} L_0\eta.L_0^2 \eta\\
    =& \ \alpha\int_{\R^2} (m_*.L_0\eta).\Big(L_0(m_*.L_0\eta) - L_0m_*.L_0\eta + 2\grad m_*:\grad L_0\eta\Big)\\
    &-\alpha\int_{\R^2}\big|\grad\Delta\eta\big|^2 \\
    =&  \ \alpha \int_{\R^2}\big|\grad (m_*.L_0\eta)\big|^2 +\alpha\int_{\R^2} (m_*.L_0\eta).\Big(2\grad m_*:\grad L_0\eta - L_0m_*.L_0\eta\Big)\\
    & -\alpha\norm{\eta}{\dot{H}^3}^2
\end{align*}

The maximal order of derivatives in the last integral above is at most $5$, from which we infer
\begin{equation*}
    \mathcal{J}_1 = -\alpha\norm{\eta}{\dot{H}^3}^2 + \alpha\norm{\grad (m_*.L_0\eta)}{L^2}^2  +\int_{\R^2}\mathcal{O}_5^2(\eta).
\end{equation*}

Finally, we compute
\begin{align*}
    \alpha\norm{\grad (m_*.L_0\eta)}{L^2}^2 &=\alpha\Big(\norm{\partial_1 (m_*.L_0\eta)}{L^2}^2 + \norm{\partial_2 (m_*.L_0\eta)}{L^2}^2\Big)\\
    &=\alpha\int_{\R^2} \Big(\big(m_*.\partial_1 L_0\eta\big)^2 + \big(m_*.\partial_2 L_0\eta\big)^2\Big) + \int_{\R^2} \mathcal{O}_5^2(\eta).
\end{align*}

Using the expression of $\Delta\eta$ in the Proof of Lemma~\ref{lemme: expression Leta} and the fact that $m_*.\partial_1 m_* = 0$, we deduce that
\begin{equation*}
    m_*.\partial_1 \Delta\eta = \partial_1(\Delta\eta.m_*) + (\Delta\eta.n_*) (m_*.\partial_1 n_*) + (\Delta\eta.p_*)( m_*.\partial_1 p_*),
\end{equation*}
that can be expressed as $\mathcal{O}_2^1(\eta) + \mathcal{O}_3^2(\eta)$.

An analogous estimate can be provided for $m_*.\partial_2 \Delta\eta$, which leads to
\begin{equation}\label{estimée J_1}
    \mathcal{J}_1 = -\alpha\norm{\eta}{\dot{H}^3}^2    +\int_{\R^2}\mathcal{O}_5^2(\eta) + \int_{\R^2}\mathcal{O}_5^3(\eta)+ \int_{\R^2} \mathcal{O} _6^4(\eta).
\end{equation}

As for the second part, it contains only lower order terms, and this reads 
\begin{equation}\label{estimée J_2}
    \mathcal{J}_2 = \int_{\R^2} \mathcal{O}_5^2(\eta).
\end{equation}

Regarding the remaining terms, man still needs to deal with each term in the right-hand side of~\eqref{expression de J}. As for the first one, we have
\begin{equation*}
    \left|-\alpha\int_{\R^2}\Big(m_*\times \big(\eta\times \heff(m_*)\big)\Big).L^2\eta\right|=\int_{\R^2} \mathcal{O}_4^2(\eta).
\end{equation*}

The terms in the second line of~\eqref{expression de J} can be shown to behave as 
\begin{equation*}
    \int_{\R^2}\mathcal{O}_4^3(\eta) + \int_{\R^2}\mathcal{O}_6^3(\eta).
\end{equation*}

As for the last line of~\eqref{expression de J}, we have 
\begin{equation*}
    \left|\alpha\int_{\R^2}\Big(\eta\times\big(\eta\times L\eta\big)\Big).L^2\eta\right|=\int_{\R^2}\mathcal{O}_6^4(\eta).
\end{equation*}

We combine~\eqref{estimée J_1},~\eqref{estimée J_2} with all remaining terms, which leads to 
\begin{align*}
    \mathcal{J} + \alpha\norm{\eta}{\dot{H}^3} ^2= \sum_{\substack{(l,k)\in \mathcal{K}}}\int_{\R^2}\mathcal{O}_k^l(\eta)
\end{align*}

with $\mathcal{K}=\big\{(4,2),(5,2),(4,3),(5,3),(6,3),(6,4)\big\}$. We invoke Lemma~\ref{lemme: controle int O^l_3 et O^l_4} to deal with all previous terms in the sum over $\mathcal{K}$ similarly to the proof of Claim~\ref{claim: controle de I} and this provides the desired estimate.
\end{proof}

We collect both estimates obtained in Claim~\ref{claim: controle de I} and~\ref{claim: controle de J + alpha norm H^3 eta}. The desired estimate comes from imposing $\varepsilon$, and then $\delta_3$ small enough so that, by Proposition~\ref{proposition: propriétés de f}, the factor multiplying $\norm{\eta}{\dot{H}^3}$ in the right-hand terms can be absorbed. There exists $C_3>0$ such that
\begin{equation*}
        \mathcal{I}+\mathcal{J} + \dfrac{\alpha}{2}\norm{\eta}{\dot{H}^3}^2\leq C_3\Big(\norm{\eta}{H^2}^2 +\norm{\eta}{H^2}^4\Big) .
    \end{equation*}

\end{proof}

\section{Control of the perturbation \texorpdfstring{$\eta$}{eta}}\label{section: control of the perturbation}

This section lays the foundation for the last bootstrap argument, which eventually yields asymptotic stability. In Subsection~\ref{section: spectral property of mathcal L}, we describe the spectrum of $\q L$, which enables us in Subsection~\ref{section: Orthogonal decomposition of the magnetization} to perturb a skyrmion while avoiding the negative spectrum of the linearized operator $\q L$.

\subsection{Spectral property of the linearized operator \texorpdfstring{$\q L$}{L} and \texorpdfstring{$L$}{L}}\label{section: spectral property of mathcal L}

\index{operator $\q L$}
\index{operator $L$}

By the translation invariance of $E$, the quadratic form $H_*$ associated with the linearized operator\footnote{According to Remark~\ref{remarque: différence entre les opérateurs linéarizés L}, the linearized operator associated with the energy is $\q L$ or $L-\Lambda_{*}$ depending on the point of view. This paragraph highlights formally that $\ker\big(L-\Lambda_{*}\big)\neq \{0\}$, but we can also refer to Subsection 2.5 in~\cite{LiMelcher} for some rigorous explanation of this fact, using geometrical features of the sphere.} about the skyrmion $m_*$ defined by $H_*(\eta_1,\eta_2)=\psLdeux{\big(L-\Lambda_{*}\big)\eta_1}{\eta_2}$ has a non-trivial kernel. Loosely speaking, this fact can be viewed by differentiating and evaluating the (shifted) Euler-Lagrange equation
\begin{equation*}
    \delta E\big(\tau_t m_*) = \Lambda (\tau_t m_*) \tau_t m_*,
\end{equation*}

with respect to time at $t=0$, where we have set $\tau_tm_* = m_*(.-tc)$ for $c\in\R^2\backslash\{0\},t>0$. It yields $\big(L-\Lambda_{*}\big)\big((c.\grad)m_*\big)=\partial_t\big[\Lambda(\tau_t m_*)\big]_{t=0}m_*$. Thus, for orthogonality reasons, we obtain that $(c.\grad)m_*\in\ker H_*$ for any $c\in\R^2$. Now, introducing the notations
\begin{align*}
    X_1:=\partial_1 m_* &= \dfrac{\sin\psi\sin\theta}{r}n_*+\cos\psi\ \theta'p_*:= \nu_1 n_* + \varpi_1 p_*\\
    X_2:=\partial_2 m_* &= -\dfrac{\cos\psi\sin\theta}{r}n_*+\sin\psi\ \theta'p_*:= \nu_2 n_* + \varpi_2 p_*,
\end{align*}
we give a rigorous description of the kernel of the linearized operator $\q L$.

\begin{lem}\label{lemme: coercivité de L}
The operator $\q L$ is a self-adjoint operator on $L^2(\R^2)\times L^2(\R^2)$ with domain $H^2(\R^2)\times H^2(\R^2)$. This operator is non-negative, $0$ is a double eigenvalue and we have
\begin{equation*}
    \ker(\q L)=\spanned{X_1,X_2}\quad\text{with }X_j=(\nu_j,\varpi_j)\text{ for }j\in\{1,2\}.
 \end{equation*}

Moreover, there exists $\lambda_0 >0$ such that for any $\eta\in H^2(\R^2)$ that decomposes as $\eta = \mu m_* + \nu n_* + \varpi p_*$ and such that $\psLdeuxLdeux{(\nu,\varpi)}{X_j}=0$ for $j\in\{1,2\}$, we have
\begin{equation}\label{coercivité de mathcal L}
    \lambda_0 \norm{(\nu,\varpi)}{H^1\times H^1}^2\leq \psLdeuxLdeux{\q L(\nu,\varpi)}{(\nu,\varpi)},
\end{equation}

and
\begin{equation}\label{coercivité H^2}
    \lambda_0 \norm{(\nu,\varpi)}{H^1}\leq \normLdeuxLdeux{\q L(\nu,\varpi)}.
\end{equation}

\end{lem}

\begin{proof}
The facts that $\q L$ is a well-defined operator on the domain $H^2(\R^2)^2$, non-negative and \eqref{coercivité de mathcal L} are proven in Proposition~4 and Proposition~5 in~\cite{LiMelcher}. These properties rely on sharp estimates made on the profile of the skyrmion, dealt with in the Appendix of the latter article. Regarding~\eqref{coercivité H^2}, we just use the Cauchy-Schwarz inequality and~\eqref{coercivité de mathcal L} to obtain
\begin{equation*}
    \lambda_0 \norm{(\nu,\varpi)}{H^1\times H^1}^2\leq \psLdeuxLdeux{\q L(\nu,\varpi)}{(\nu,\varpi)}\leq  \normLdeuxLdeux{\q L(\nu,\varpi)}\norm{(\nu,\varpi)}{H^1},
\end{equation*}
and we conclude by simplification.

Arguing as in~\cite{LiMelcher}, we can show that the kernel of $\q L$ is made of the tangent plane to $\S^2$ at $m_0$, i.e. the two-dimensional plane in $\R^3$ given by $\spanned{\partial_1 m_0,\partial_2 m_0}$ (since $\partial_1 m_0$ and $\partial_2 m_0$ are linearly independent). For sake of completeness, we compute
\begin{align*}
   \q L_\nu\nu_1 + S \varpi_1 &= -\Delta\nu_1 + \big(N_\nu-\Lambda_{*}\big)\nu_1 - P_{\nu,\varpi}\partial_\psi\varpi_1\\
   &=-\sin\psi\left(\dfrac{\theta''\cos\theta}{r}-\dfrac{(\theta')^2\sin\theta}{r}-\dfrac{\theta'\cos\theta}{r^2}\right) \\
   &\quad\quad\quad+\left(\dfrac{1}{r^2} + h - \Lambda_{*}\right)+ \dfrac{2}{r}\left(\sin\theta - \dfrac{\cos\theta}{r}\right)\sin\psi\ \theta'.
\end{align*}

Plugging the expressions~\eqref{équation sur theta} and~\eqref{expression multiplicateur de Lagrange lambda m_*} into the expression above yields the cancelation of the latter. Second, we have
\begin{align*}
    S^*\nu_1 + \q L_\varpi \varpi_1 =& P_{\nu,\varpi}\partial_\psi\nu_1  - \Delta\varpi_1 + \big(P_\varpi - \Lambda_{*}\big)\varpi_1\\
    =& \dfrac{2}{r}\left(\sin\theta - \dfrac{\cos\theta}{r}\right)\dfrac{\cos\psi\sin\theta}{r} - \cos\psi\left(\theta''' + \dfrac{\theta''}{r}-\dfrac{\theta'}{r^2}\right)\\
    &\quad\quad\quad + \ \theta'\cos\psi\left(-\dfrac{4\sin\theta\cos\theta}{r}+\dfrac{\cos^2\theta - \sin^2\theta}{r^2}+h\cos\theta\right).
\end{align*}
By differentiating~\eqref{équation sur theta} and plugging the expression of $\theta'''$ in the expression above, the formula simplifies and leads to $\q L(\nu_1,\varpi_1)=0$.
\end{proof}

\begin{rem} \label{rem:L_H^2}
One should emphasize that, somewhat surprisingly, there is no $H^2$ control of $\q L$: there is no constant $c>0$ such that for all $\nu,\varpi \in H^2$,
\[ \| \q L(\nu, \varpi) \|_{L^2} \ge c \| (\nu,\varpi) \|_{H^2}. \]
Indeed, inspecting the definition of $\q L_\nu$ and $\q L_\varpi$, we see that $\Lambda_*$ is bounded but $N_\nu$ and $P_\varpi$ are equivvalent to $1/r^2$ as $r \to 0$: so these operators (and thus $\q L$) are $H^1$ perturbations of $- \Delta + \frac{1}{r^2}$. Now, the inverse square potential has critical singularity in space dimension $2$, and comes here with critical size $1$. 

Nevertheless, we will be able to recover an $H^2$ control using $L\eta$. Its coercivity over $H^2$ is the content of the next lemma.
\end{rem}

\index{operator $L$}

\begin{lem}\label{lemme: Uniform coercivity of H_m_0(eta) for large h}
There exists $\lambda_1>0$, only depending on $h>1$, such that for any $\eta\in H^2(\R^2)$, we have
\begin{equation}\label{eq:coercivity}
\lambda_1 \|\eta\|_{H^2} \;\le\;  \, \| L\eta\|_{L^2}\leq \frac{1}{\lambda_1} \norm{\eta}{H^2} .
\end{equation}  
\end{lem}

\begin{proof}[Proof of Lemma~\ref{lemme: Uniform coercivity of H_m_0(eta) for large h}]
We work in Fourier variables. As 
\[ \widehat L\eta (\xi) = (|\xi|^2 + h)\widehat{\eta}(\xi)+ 2i\xi\times \widehat{\eta}(\xi) = M(\xi) \widehat{\eta}(\xi) ,\]
with
\begin{equation*}
    M(\xi)=\begin{pmatrix}
        |\xi|^2+h & 0 & 2i\xi_2 \\
        0 & |\xi|^2+h & -2i\xi_1\\
        -2i\xi_2 & 2i\xi_1 & |\xi|^2+h
    \end{pmatrix},
\end{equation*}
it suffices to prove that $|M(\xi)\widehat{\eta}(\xi)|^2 \gtrsim \langle\xi\rangle^2|\widehat{\eta}(\xi)|^2$.

We compute that $\sigma(M)=\{|\xi|^2+h,|\xi|^2+h\pm 2|\xi|\}$ and then the smallest eigenvalue is then $\lambda_{\min}(|\xi|)=|\xi|^2 -2|\xi| + h$. Thus, by standard spectral argument, and a variational study of the real function $x\mapsto \frac{x^2 - 2x + h}{x^2 + 1}$, we obtain for $h >1$,
\begin{equation*}
    (h-1)|\widehat{\eta}(\xi)|^2\leq\lambda_{\min}(|\xi|)^2|\widehat{\eta}(\xi)|^2\leq |M(\xi)\widehat{\eta}(\xi)|^2.
\end{equation*}
This gives the bound from below. The bound from above is clear.
\end{proof}

\subsection{Orthogonal decomposition of the magnetization}\label{section: Orthogonal decomposition of the magnetization}

The non-trivial kernel for the linearized operator $\q L$ defined in the previous subsection is an actual obstruction to obtain coercive properties. We modulate the solution in such a way that the solution $m(t)$ will decompose into the skyrmion $m_*$ plus an extra perturbation, orthogonal to the kernel of the non-negative operator $\q L$, then avoiding the unstable spectrum of the latter operator. In this matter, we define a modulated tubular neighborhood of a vector $m\in H^s_*$, and we will check that, the magnetization $m$ decomposes as previously mentioned, as long as it remains in such a tube of sufficient small depth. We define the set\index{$\mathcal{U}(\delta)$}
\begin{equation}\label{voisinage tubulaire autour de m_* avec translation}
    \mathcal{U}(\delta) = \left\{m\in H^2_* \mid \inf_{x\in\R^2}N_f\big(m(.-x),m_*\big)<\delta\right\}\index{$N_f(\eta)$}.
\end{equation}

\begin{lem}\label{lemme: décomposition orthogonale}
    There exist $0<\delta_4\leq \delta_3,C_ 4>0$ and a function $x=(x_1,x_2)\in C^1\big(\mathcal{U}(\delta_4),\R^2\big)$ such that the following holds. For any $0<\delta\leq \delta_4$ and any $m\in \mathcal{U}(\delta)$, the function defined by 
\begin{equation*}
    \eta:=\eta(m)=m\big(.-x(m)\big) - m_*,
\end{equation*}\index{orthogonal decomposition}
satisfies the orthogonality conditions
\begin{equation*}
    \psLdeux{\eta}{\partial_1 m_*} = \psLdeux{\eta}{\partial_2 m _ *}=0, 
\end{equation*}
    and\index{$N_f(\eta)$}
    \begin{equation} \label{est:N_f_mod}
        N_f\big(\eta(m)\big) + \normR{x(m)}\leq C_4\delta,
    \end{equation}
    as well as the pointwise identity on $\R^2$, that is
    \begin{equation*}
        |\eta|^2 + 2 \eta.m_* = 0.
    \end{equation*}
    
\end{lem}

\begin{proof}
    First of all, consider the mapping defined for any $(m,x)\in  H^1(\R^2,\R^3)\times \R^2$ by
    \begin{equation}
        X(m,x)=\begin{pmatrix}
            \psLdeux{m(.-x)-m_*}{\partial_1 m_*}\\
            \psLdeux{m(.-x)-m_*}{\partial_2 m_*}
        \end{pmatrix}.
    \end{equation}

We first restrict $X$ to $\mathcal{M}\times \R^2$ with $\mathcal{M
}:=\big\{m\in H^1_\loc(\R^2,\R^3)\big|\exists\eta\in H^1(\R^2,\R^3),m=m_* + \eta\text{ and }|m|^2 = 1\big\}$. Similarly to the proof of Lemma~4.2 in~\cite{CotIgnat}, we can check that the restriction $X: \mathcal{M}\times \R^2\rightarrow\R^2$ defines a $C^\ii$ mapping on a submanifold of $H^1(\R^2,\R^3)\times \R^2$. We have $X(m_*,0)=0$ and \begin{equation}\label{déterminant de partial_2 X non nul}
    \partial_2 X(m_*,0)=-\begin{pmatrix}
        \normLdeux{\partial_1 m_*}^2 & \psLdeux{\partial_1m_*}{\partial_2 m_*}\\
        \psLdeux{\partial_1 m_*}{\partial_2 m_*} & \normLdeux{\partial_2 m_*}^2
    \end{pmatrix}\in GL_2(\R).
\end{equation}

Indeed, $\det\big(\partial_2 X(m_*,0)\big)=\normLdeux{\partial_1 m_*}^2 \normLdeux{\partial_2 m_*}^2 -\psLdeux{\partial_1m_*}{\partial_2 m_*}^2 $ and if the determinant was zero, this would correspond to the equality case in the corresponding Cauchy-Schwarz inequality. Thus, $\partial_1 m_*$ and $\partial_2 m_* $ should be linearly dependent, which contradicts the fact that the skyrmion number is not zero, i.e. that $Q(m_*)\neq 0$\index{skyrmion number}. We now invoke a modified implicit function argument, that is reminiscent from~\cite{BetGrSm1} (see the appendix therein). The first deduction of~\eqref{déterminant de partial_2 X non nul} is the existence of an $H^1_*$ neighborhood $V$ of $m_*$ and a function $\overline{x}:V\rightarrow\R^2$ such that $X\big(m,\overline{x}(m))=0$ for any $m\in V$. We then prove that the assumptions of Proposition A.1 in~\cite{BetGrSm1} are satisfied with $E\times F=\mathcal{M}\times \R^2,\phi=X$.

In fact, we have $\big|\partial_2 X(0,m_*)\big|\leq \max\big(\normLdeux{\partial_1 m_*},\normLdeux{\partial_1 m_*}\big)$ and we check that for any ball $\mathcal{B}\subset H^1(\R^2,\R^3)$, $\grad^2 X$ is bounded on $\mathcal{B}\times \R^2$. Proposition A.1 in~\cite{BetGrSm1} applies and there exist $\beta,\Lambda$ such that for any $x\in\R^2$, we can define a function $\gamma_x\in C^1\big(\mathcal{B}(m_*,\beta),\R^2\big)$ such that for any $m\in \mathcal{B}(m_*,\beta)$, $\gamma_x(m)$ is the only solution in $\mathcal{B}(x,\Lambda \beta)$ to $X(m,\gamma_x(m))=0$, with $\gamma_x(m_*)=0$ and $\gamma_x$ is Lipschitz continuous on $\mathcal{B}(m_*,\beta)$.

Assume that $\delta_4$ is small enough so that $2\delta_4 <\beta^2$. If $m\in \mathcal{U}(\delta)$, and by Proposition~\ref{proposition: propriétés de f}, there exists $x\in \R^2$ such that $\| m(.-x),m_* \|_{H^1} \leq N_f ( m(.-x) - m_*) <2\delta$. Thus $m(.-x)\in\mathcal{B}(m_*,2\delta)$. Now set $x(m):=\gamma_x(m)$. By the Lipschitz continuity property, we have 
\begin{equation*}\label{controle x(m) dans fonction implicite}
    |x(m)|=|x(m)-x(m_*)|\leq \Lambda \| m- m_* \|_{H^1} \leq \Lambda B_* f\big(m-m_*)\leq 2 \Lambda B_* \delta. 
\end{equation*}
Observe that $\| (m_*(.+y_1) - m_*(.+y_2) \|_{H^2} \leq  K|y_1-y_2|$ where $K$ depends only on $\norm{m_*}{\mathcal{Z}^k}$\index{$\mathcal{Z}^{s}(\R^d)$}. In particular, $\| m_*(.+x(m)) -m_* \|_{H^1} \leq K|x(m)|\leq K\Lambda \| m - m_* \|_{H^1}$, hence by Proposition~\ref{proposition: propriétés de f}, there exists $C_4>0$ such that
\begin{align*}
     N_f\big(\eta(m)\big)
    &\le f\Big( \| m -m_* \|_{H^1} + \| m_*(.+x(m) -m_* \|_{H^1} , \\
    & \qquad \qquad \| m - m_* \|_{H^2} + \| m_*(.+x(m)\big) - m_* \|_{H^2}) \le C_4\delta. \qedhere
\end{align*}

\end{proof}

\section{Proof of the main results}\label{section: proof of the main results}

We now prove Theorems~\ref{théorème: LWP of the Cauchy problem} and~\ref{théorème: main stabilité asymptotique skyrmion}. First, we handle the local well-posedness of~\eqref{llg}. This shall enable us to differentiate several quantities and then use the dissipative features of~\eqref{llg}. Then, a proper arrangement of these previous quantities will provide enough control to close a bootstrap argument, and eventually conclude as for the stability result. 

\bigskip

First recall that we stated the Cauchy problem for~\eqref{llg} in a more general framework, for dimensions $d=2$ and $d=3$, in Theorem \ref{théorème: local well-posedness} and \ref{théorème: LWP of the Cauchy problem}. The proofs are done in Appendix~\ref{appendix: LWP}.

Then in Subsection~\ref{section: Bootstrap argument for a smooth initial data} we give bound on the dynamics for the orthogonal decomposition of the magnetization $m(t)$, and then we Theorem~\ref{théorème: main stabilité asymptotique skyrmion}, for smooth initial data, allowing us to use the estimates of Section~\ref{section: control of the perturbation}. This is the core of the proof. At last in Subsection~\ref{End of the proof for an initial data in}, we recover the stability in the original setting by a density argument.

%
%
%
%
%
%
%
%
%

\subsection{Bootstrap argument for a smooth initial data}\label{section: Bootstrap argument for a smooth initial data}

As we rely on dissipation identities  \eqref{dérivée de E(m(t)) non rigoureuse} and \eqref{dérivée ordre deux non rigoureuse}, we need an extra regularity on $m_0$, namely $m_0\in H^3_*$ to perform the estimates; this condition will be relaxed in the last step by an approximation argument  in Subsection~\ref{End of the proof for an initial data in}.


First let us start with a dynamical bound on the translation parameter, for a magnetization solution of \eqref{llg}, and which admits a decomposition as in Lemma \ref{lemme: décomposition orthogonale} on a whole interval of time.

\begin{lem}[Dynamical modulation of the magnetization]\label{lemme: décomposition orthogonale dynamique}
    There exist $0<\delta_6 \leq\delta_4$ and $C_6>0$ such that the following holds. 
    Let $m \in \q C([0,T], H^2_*(\m R^2))$ be a solution to \eqref{llg} on $[0,T]$, such that
    \[ \forall t \in[0,T], \quad m(t) \in \mathcal{U}(\delta), \]
    so that by Lemma \ref{lemme: décomposition orthogonale}, it admits a decomposition 
    \[ m(t,x) = m_*(x-x(t)) + \eta(t) \]
    where  $x\in C^1\big([0,T],\R^2\big)$, $\eta \in \q C([0,T], H^2(\m R^2))$ (together with the orthogonality properties and bound of Lemma \ref{lemme: décomposition orthogonale}. Then, for any $t\in [0,T]$, there hold
    \begin{equation}
        |\dot{x}(t)|\leq C_6\normHun{\eta(t)},
    \end{equation}
\end{lem}

\begin{proof}
In view of the smoothness of the initial data, $\eta$ can be differentiated with respect to time and takes values in $L^2(\R^2)$. For similar reasons, $x$ is actually in $C^1\big([0,T],\R^2\big)$. Therefore, using the Euler-Lagrange equation~\eqref{euler lagrange equation point critique}, we obtain
\begin{align}
    \partial_t\eta(t) =&\ \dot{x}(t).\grad m_*\big(.-x(t)\big) +  \eta(t)\times\heff\big(m(t)\big)-m_*\times L\eta(t)\label{llg sur eta}\\
    & \quad\quad\quad \quad -\alpha m(t)\times\Big(\eta(t)\times\heff\big(m(t)\big)-m_*\times L\eta(t)\Big)\notag.
\end{align}

Differentiating the orthogonality conditions $\psLdeux{\eta(t)}{\partial_j m_*\big(.-x(t)\big)}=0$ for $j\in\{1,2\}$ with respect to $t$ and plugging~\eqref{llg sur eta} into the latter writes in terms of $2\times 2$ matrices $M=(m_{i,j})_{1\leq i,j\leq 2}, H(t)$ and a two-dimensional vector $\Phi(t)=\big(\Phi_1(t),\Phi_2(t)\big)$, and reads
\begin{equation}
    \widetilde{M}(t)\begin{pmatrix}
        \dot{x}_1(t)\\
        \dot{x}_2(t)
    \end{pmatrix}=\Phi(t),
\end{equation}

with $\widetilde{M}(t)=M + H(t)$ defined by
\begin{equation*}
    M=\begin{pmatrix}
        \normLdeux{\partial_1 m_*}^2 & \psLdeux{\partial_1m_*}{\partial_2 m_*}\\
        \psLdeux{\partial_1m_*}{\partial_2 m_*} & \normLdeux{\partial_2 m_*}^2
    \end{pmatrix},
\end{equation*}
and
\begin{equation*}
        H(t)=\begin{pmatrix}
        \psLdeux{\partial_1\eta(t)}{\partial_1 m_*\big(t,.-x(t)\big)} & \psLdeux{\partial_2\eta(t)}{\partial_1 m_*\big(t,.-x(t)\big)}\\
        \psLdeux{\partial_1\eta(t)}{\partial_2 m_*\big(t,.-x(t)\big)} & \psLdeux{\partial_2\eta(t)}{\partial_2 m_*\big(t,.-x(t)\big)}
    \end{pmatrix},
\end{equation*}
and 
\begin{align*}
    \Phi_j(t)&=\psLdeux{\eta(t)\times\heff\big(m(t)\big)-m_*\times L\eta(t)}{\partial_j m_*}\\
    &-\alpha \psLdeux{ m(t)\times\Big(\eta(t)\times\heff\big(m(t)\big)-m_*\times L\eta(t)\Big)}{\partial_j m_*\big(t,.-x(t)\big)}.
\end{align*}

Arguing as in the proof of Lemma~\ref{lemme: décomposition orthogonale}, we verify that the determinant of $M$ is not zero and we write $\widetilde{M}(t)= M\big(I + M^{-1}H(t)\big)$, and we can check that there exists a constant $C>0$ such that $|\Phi_1(t)|+|\Phi_2(t)|+|H(t)|\leq C\normHun{\eta(t)}$. From Lemma~\ref{lemme: décomposition orthogonale}, we also have $N_f\big(\eta(t)\big)\leq C_4\delta$. By continuity of the determinant, there exists $\delta_6$ small enough so that for $\delta\leq \delta_6,|M^{-1} H(t)|<1$ uniformly on $[0,T]$ and $\widetilde{M}(t)$ is invertible and $|\widetilde{M}(t)^{-1}|$ is bounded by a constant independent of $t$. To conclude, there exists a constant $C_6>0$ such that for any $t\in [0,T],|\dot{x}(t)|\leq C_6\normHun{\eta(t)}.$  
\end{proof}

Here, we have $\gm:=m_*\in\mathcal{Z}^3(\R^2)$\index{$\mathcal{Z}^{s}(\R^d)$} which satisfies~\eqref{Z_0 in H^s-1} with $s=3$, due to Proposition~\ref{proposition: propriétés du skyrmion}. By Theorem~\ref{théorème: local well-posedness}, we can define $m=\gm+w$ the solution associated with the initial data $m_0$, defined on the maximal interval of existence $[0,T_{\max})$. To suit the assumptions of the main theorem, assume that \index{$N_f(\eta)$}
\begin{equation}\label{initial perturbation d_2}
    N_f(m_0-m_*):=\delta< \delta_*,
\end{equation}
with some $\delta_*$ which is fixed according to several conditions, among which~\eqref{C_6Mdelta_* leq delta_3},\\~\eqref{choix delta_*},~\eqref{choix delta_* deyx} and~\eqref{condition sur delta_* trois}, which are highlighted below in the proof of the next proposition. In addition, we further anticipate by imposing that $\delta_*$ is small enough so that $ M\delta_*<\delta_6$, where $M$ is a fixed large constant set in the following proposition (see~\eqref{valeur de M}) and most importantly independent upon $\delta_*$. As a matter of fact, this special choice of $\delta_*$ enables us to use the orthogonal decomposition of $m(t)$ and the estimates in Lemma~\ref{lemme: décomposition orthogonale dynamique}, that are crucial in the sequel.

\begin{proof}[Proof of Theorem~\ref{théorème: main stabilité asymptotique skyrmion} with $m_0\in H^3_*$.]

Assume $M\geq 2C_4$, as imposed in~\eqref{valeur de M}. By~\eqref{initial perturbation d_2}, $m_0\in \mathcal{U}(\delta_*)$, thus by Lemma~\ref{lemme: décomposition orthogonale}, 
\begin{equation}\label{condition initiale C_4delta_*}    
N_f\big(\eta(0)\big)\leq C_4\delta_*\leq \dfrac{M\delta_*}{2}.
\end{equation}

By~\eqref{condition initiale C_4delta_*}, we can define \index{$\mathcal{U}(\delta)$}
\begin{equation}
    T_*=\sup\big\{0<T< T_{\max} \mid \forall t \in [0,T],m(t)\in\mathcal{U}(M\delta_*)\big\}>0.
\end{equation}

Since by assumption $M\delta_* <\delta_6$, Lemma~\ref{lemme: décomposition orthogonale dynamique} applies on $[0,T]$ and we have \begin{equation}
 \forall t \in [0,T], \quad   f\Big(\norm{\eta(t)}{H^1},\norm{\eta(t)}{H^2}\Big)\leq C_4 M\delta_*.
\end{equation}

One of the choices of the value of $\delta_*$ is made according the following condition. Given $M$, we impose indeed that \begin{equation}\label{C_6Mdelta_* leq delta_3}
    C_4 M\delta_*\leq \delta_3.
\end{equation}

For $t\in [0,T]$, we set the time-dependent functions $p_j(t)=\norm{\eta(t)}{H^j}^2$ for $j\in\{1,2,3\},\ q(t)=E\big(m(t)\big)-E(m_*)$ and $s(t)=\norm{L\eta(t)}{L^2}^2$, and we gather the dissipative estimates of the previous sections. Ultimately, we show how they articulate in order to obtain the exponential decay of the $H^2$-norm of the perturbation $\eta(t):=\mu(t) m_* + \nu(t)n_* + \varpi(t) p_*$.

Combining Lemma~\ref{lemme: controle estimation E(m_*+eta)-E(m_*)} with Lemma~\ref{lemme: coercivité de L}, we have for any $t\in [0,T]$,
\begin{equation*}
    -C_1\delta_1p_1(t) +\dfrac{\lambda_0}{2}\norm{\big(\nu(t),\varpi(t)\big)}{H^1\times H^1}^2 \leq q(t)\leq C_1p_1(t) .
\end{equation*}

Moreover, using Lemma~\ref{lemme: equivalence norm H^1 eta et (pi,nu)}, and up to reducing the value of $\delta_*$ so that 
\begin{equation}\label{choix delta_*}
    \dfrac{\lambda_0}{4C_0}\leq\dfrac{\lambda_0}{2C_0}-C_1\delta_*,
\end{equation}
we obtain a constant $\lambda_1>0$ such that for any $t\in [0,T]$,
\begin{equation}\label{équivalence entre p_1 et q}
    \dfrac{1}{\lambda_1}p_1(t)\leq q(t)\leq \lambda_1 p_1(t).
\end{equation}

Second, combining Lemma~\ref{lemme: equivalence norm H^1 eta et (pi,nu)},~\ref{lemme: controle perturbation niveau H^1} and~\ref{lemme: coercivité de L}, there exists $\lambda_2>0$ such that
\begin{equation*}
    \dot{q}(t)  \leq -\alpha\norm{\q L\big(\nu(t),\varpi(t)\big)}{L^2\times L^2}^2 + C_2\delta_2\normHun{\eta(t)
    }^2\leq \left(C_2\delta_2-\dfrac{\alpha\lambda_0^2}{C_0}\right)\normHun{\eta(t)}^2.
\end{equation*}

Thus, imposing $\delta_*$ small enough so that 
\begin{equation}\label{choix delta_* deyx}
    C_2\delta_2-\frac{\alpha\lambda_0^2}{C_0}\leq -\frac{\alpha\lambda_0^2}{2C_0}=:-\lambda_2<0,
\end{equation}
we infer from~\eqref{équivalence entre p_1 et q} that for any $t\in [0,T]$,
\begin{equation}
    \dot{q}(t)\leq -\dfrac{\lambda_2}{\lambda_1}q(t), 
\end{equation}

hence, by~\eqref{équivalence entre p_1 et q},
\begin{equation}\label{p_1(t) décroit expo}
    p_1(t)\leq \lambda_1^2 \ p_1(0)e^{-\frac{\lambda_2}{\lambda_1}t}.
\end{equation}

From the initial control~\eqref{condition initiale C_4delta_*} and Proposition~\ref{proposition: propriétés de f}, we have for any $t\in [0,T]$,
\begin{align}\label{p_1(t) décroit expo avec delta_*}
    p_1(t)\leq  \lambda_1^2 B_*^2 f\Big(\sqrt{p_1(0)},\sqrt{p_2(0)}\Big)^2e^{-\frac{\lambda_2}{\lambda_1}t}\leq\lambda_1^2 B_*^2C_4^2\delta _ *^2e^{-\frac{\lambda_2}{\lambda_1}t}.
\end{align}

From Lemma~\ref{lemme: Uniform coercivity of H_m_0(eta) for large h}, we can derive the fact that there exists $\lambda_3>0$ such that for any $t\in [0,T]$, we have \begin{equation}\label{équivalence s et p_2}
    \dfrac{1}{\lambda_3}p_2(t)\leq s(t)\leq \lambda_3 p_2(t),
\end{equation}

and one can rigorously use the estimate~\eqref{dérivée ordre deux non rigoureuse}, which reads
\begin{equation*}
    \dot{s}(t)\leq Diss_2\big(m(t)\big)+C_H|\dot{x}(t)|\norm{m_*}{Z^3}\normLdeux{L\eta(t)}.
\end{equation*}

From Lemma~\ref{lemme: controle au niveau H^2 perturbation}, Lemma~\ref{lemme: Uniform coercivity of H_m_0(eta) for large h} and Lemma~\ref{lemme: décomposition orthogonale dynamique}, we infer
\begin{equation*}
    \dot{s}(t) + \dfrac{\alpha}{2}\norm{\eta(t)}{\dot{H}^3}^2\leq C_3\big(p_2(t) + p_2(t)^2\big) + \frac{2C_6C_H}{\lambda_1} \norm{m_*}{Z^3}\sqrt{p_1(t)p_2(t)}.
\end{equation*}

Writing $\widetilde{C}_3 = C_3 + \frac{2C_6C_H}{\lambda_1} \norm{m_*}{Z^3} \geq C_3$, and interpolating the $H^2$-norm between $H^1$ and $H^3$, we have 
\begin{align*}
    \dot{s}(t) + \dfrac{\alpha}{2}p_3(t) &\leq \left( \widetilde{C}_3 + \dfrac{\alpha}{2}\right)p_2(t) + \widetilde{C}_3p_2(t)^2\\
    &\leq \left( \widetilde{C}_3 + \dfrac{\alpha}{2}\right)\sqrt{p_1(t)p_3(t)}+\widetilde{C}_3p_1(t)p_3(t),
\end{align*}

Hence, for any $\varepsilon >0 $, by Young's inequality $2|ab|\leq \varepsilon a^2 + \frac{b^2}{\varepsilon}$, there exists a constant $C_\varepsilon >0$ such that
\begin{equation}
    \dot{s}(t)\leq  \Big(\dfrac{(2\widetilde{C}_3+\alpha)\varepsilon}{4}+\widetilde{C}_3p_1(t) - \dfrac{\alpha}{2}\Big) p_3(t) + C_\varepsilon p_1(t).
\end{equation}

We set $\varepsilon=\frac{\alpha}{2(2\widetilde{C}_3+\alpha)}$ and 
\begin{equation}\label{condition sur delta_* trois}
    \delta_*\leq \sqrt{\dfrac{\alpha}{8\lambda_1^ 2B_*^2\widetilde{C}_3C_4^2}},
\end{equation}

so that, taking~\eqref{équivalence s et p_2},~\eqref{p_1(t) décroit expo} and~\eqref{p_1(t) décroit expo avec delta_*} into account, we obtain
\begin{equation} \label{est:s_dot_final}
    \dot{s}(t)\leq -\dfrac{\alpha}{4\lambda_3}s(t) + C_\varepsilon \lambda_1^2 p_1(0)e^{-\frac{\lambda_2}{\lambda_1}t}
\end{equation}

without loss of generality, we can impose that $\lambda_3$, defined in~\eqref{équivalence s et p_2}, is even larger so that $\widetilde{\lambda}:=\frac{\alpha}{4\lambda_3}<\frac{\lambda_2}{2\lambda_1}$. Then, using~\eqref{p_1(t) décroit expo avec delta_*} and~\eqref{équivalence s et p_2}, we have
\begin{equation*}
    \dot{s}(t)\leq -\widetilde{\lambda}\ s(t) + C_\varepsilon \lambda_1^2 p_1(0) e^{-\frac{\lambda_2}{2\lambda_1}t}.
\end{equation*}

Therefore, we have 
\begin{equation}
    s(t)\leq s(0)e^{-\widetilde{\lambda}t} + \dfrac{C_\varepsilon \lambda_1^2 p_1(0)}{\frac{\lambda_2}{2\lambda_1}-\widetilde{\lambda}}\Big(e^{-\widetilde{\lambda}t} -e^{-\frac{\lambda_2}{2\lambda_1}t} \Big).
\end{equation}

Thus, by~\eqref{équivalence s et p_2}, there exists $\lambda_4>0$ and $\tau_1>0$ such that
\begin{equation}\label{p_2(t) décroit expo}
    p_2(t)\leq \lambda_3\Big(\lambda_3 p_2(0) + C_\varepsilon \lambda_1^2 p_1(0)\Big)e^{-\tau_1t},
\end{equation}

and also by~\eqref{p_1(t) décroit expo},~\eqref{p_2(t) décroit expo}, Proposition~\ref{proposition: propriétés de f}, and the standard inequality $\sqrt{a+b }\leq \sqrt{a}+\sqrt{b}$,
\begin{align*}
    f\Big(\sqrt{p_1(t)},\sqrt{p_2(t)}\Big)&\leq f\Big(\lambda_1\sqrt{p_1(0)},\lambda_3\sqrt{p_2(0)}+\sqrt{C_\varepsilon\lambda_3}\lambda_1\sqrt{p_1(0)}\Big)\\
    &\leq\max\big(\lambda_1,\lambda_3 + \sqrt{C_\varepsilon\lambda_3}\lambda_1\big)f\Big(\sqrt{p_1(0)},\sqrt{p_2(0)}\Big)
\end{align*}

then, for any $t\in [0,T]$, \begin{equation*}
    N_f\big(\eta(t)\big)=f\Big(\norm{\eta(t)}{H^1},\norm{\eta(t)}{H^2}\Big)\leq \dfrac{M\delta_*}{2},
\end{equation*}

where we have set
\begin{equation}\label{valeur de M}
    M=2\max\Big(C_4,\max\big(\lambda_1,\lambda_3 + \sqrt{C_\varepsilon\lambda_3}\lambda_1\big)\Big).
\end{equation}

Thus, the previous estimate leads to the fact that $m(t)\in \mathcal{U}(\frac{3M\delta_*}{4})$. By a continuity argument, we are forced to have $T_*=T_{\max}$ and the previous estimates are now true for any $t\in [0,T_{\max})$. In particular,~\eqref{p_2(t) décroit expo} holds for any $t\in [0,T_{\max})$ and by the blow-up criterion in Theorem~\ref{théorème: local well-posedness}, we necessarily have $T_{\max}=+\ii$. This provides a constant $C_*>0$ such that for any $t\geq 0$, 
\begin{equation*}
    \normHun{\eta(t)}\leq C_*\delta_* e^{-\tau_*t},
\end{equation*}
by~\eqref{p_1(t) décroit expo avec delta_*} and Proposition~\ref{proposition: propriétés de f}. Furthermore, by~\eqref{p_2(t) décroit expo}, and up to increasing the value of $C_*$,
\begin{equation*}
    \norm{\eta(t)}{H^2}\leq C_* \Big(\| m_0 - m_* \|_{H^2} + \delta_*\Big)e^{-\tau_*t},
\end{equation*}

and by Lemma~\ref{lemme: décomposition orthogonale dynamique},
\begin{equation}\label{dot x (t) + eta (t)}
    |\dot{x}(t)|\leq  C_*\| m_0 - m_* \|_{H^1} e^{-\tau_*t}.
\end{equation}

To finish, we evaluate the dynamics of the modulation parameter $x$. Taking any sequence $(t_n)$ tending to $+\ii$, we deduce by~\eqref{dot x (t) + eta (t)} that $\big(x(t_n)\big)$ is a Cauchy sequence and then converges to some $x_\ii\in\R$. If, by contradiction, $\big(x(s_n)\big)$ should converge to a different limit for a another sequence $(s_n)$ (also tending to $+\ii$), then passing to the limit $n\rightarrow +\ii$ in $\big|x(t_n)-x(s_n)\big|\leq\int_{s_n}^{t_n}\big|\dot{x}(r)\big|dr\lesssim \frac{C_6 \delta_*}{\tau_*}(e^{-\tau_*s_n}-e^{-\tau_*t_n})$ would lead to a contradiction. Therefore, $x(t)\rightarrow x_\ii$ as $t\rightarrow +\ii$, and up to taking a larger constant $C_*$, we have $|x(t)-x_\ii|\leq \int_t^\ii|\dot{x}(s)|ds\leq C_*\| m_0 - m_* \|_{H^1}e^{-\tau_* t}$.

\end{proof}

\subsection{Proof of Theorem~\ref{théorème: main stabilité asymptotique skyrmion}}\label{End of the proof for an initial data in}

We now finish the proof of Theorem~\ref{théorème: main stabilité asymptotique skyrmion}. Take a initial data $m_0\in H^2_*$ and $m$ the corresponding solution on $[0,T_{\max})$. By density, one can find a regularizing sequence $(m_0^n)_n \in (H^3_*)^\N$ tending to $m_0$ in $H^2_*$. The results of subsection~\ref{section: Bootstrap argument for a smooth initial data} apply to $m_0^n$. By Theorem~\ref{théorème: local well-posedness}, we can drop the first terms of the sequence so that for any $0<T<T_{\max}$, $[0,T]$ is contained in the maximal interval $[0,T_{\max}^n)$ of existence of the solution $m^n$ to~\eqref{llg} associated with the initial data $m^n_0$. In addition, we have $m^n-m_*\in C^0\big([0,T^n_{\max}),H^3(\R^2)\big)$ and for any $T\in (0,T_{\max})$,
\begin{equation*}
    m^n - m_*\underset{n\rightarrow +\ii}{\longrightarrow}m-m_*\quad\text{in }C^0([0,T],H^2(\R^2)\big).
\end{equation*}

By~\eqref{1ere hypothese f} and~\eqref{2eme hypothese f}, $N_f(m_0^n-m_0)\rightarrow 0$ as $n\rightarrow +\ii$, then, up to dropping even more terms in the sequence $(m_0^n)$, we can assume that $N_f(m_0^n-m_*)\leq 2B_*N_f(m_0-m_*)\leq 2B_*\delta_*$. Thus, $m_0^n\in H^3_*$ complies with the framework of the previous subsections, hence we have the orthogonal decomposition $m^n\big(t,.-x^n(t)\big)= m_* + \eta^n(t)$ and there exists $x_n^\ii$, such that for any $t\in [0,T_{\max})$,
\begin{equation}\label{estimée dot x ^n leq normHun eta^n}
    |\dot{x}^n(t)|\leq C_*\normHun{\eta^n(t)}
\end{equation}
as well as
\begin{align}
     \norm{\eta^n(t)}{H^2} & \le C_* \Big( \| m_0^n - m_* \|_{H^2} +\delta_*\Big) e^{-\tau_* t} \nonumber \\
    &\leq 2B_*C_* \Big(\| m_0 - m_* \|_{H^2}+\delta_*\Big)e^{-\tau_* t}. \label{estimée x^n - c^ii_n + normHdeux eta^n}
\end{align}

Based on previous estimates, the blow-up criterion in Theorem~\ref{théorème: local well-posedness} implies that $T_{\max}=+\ii$ and for any $t\in [0,\ii)$. From the Arzela-Ascoli theorem, and up to extracting a subsequence, we obtain a limit function $x\in C^0_{\loc}(\R_+,\R^2)$ and $x_\ii\in\R^2$ such that
\begin{equation*}
    x^n\underset{n\rightarrow +\ii}{\longrightarrow}x\quad\text{in }\ C^0_{\loc}(\R_+,\R^2) \quad\text{and}\quad x^\ii_n\underset{n\rightarrow +\ii}{\longrightarrow}x_\ii.
\end{equation*}

Moreover, up to taking a further subsequence, $(\eta^n)$ converges towards a limit $\eta$ in $C^0_\loc\big(\R_+,H^2(\R^2)\big)$. Now passing to the limit $n\rightarrow +\ii$ in the orthogonal decomposition of $m^n$ and then in~\eqref{estimée x^n - c^ii_n + normHdeux eta^n} yields $m\big(t,.-x(t)\big) = m_* + \eta(t)$ and the desired estimates.

\qed

\appendix

\section{Symmetries}\label{appendix: symétries}

Both Dirichlet and Zeeman energies are invariant under spatial action of the orthogonal group $SO(2)$ and the combined reflections of the type~\eqref{combined reflection}, namely for a mapping $m:\R^2\rightarrow\S^2$, one has for any $R_1,R_2\in SO(2)\cup\{\pm F\},\ E_Z\big(\widetilde{R}_1m(R_2.)\big)=E_Z(m)$, and $E_{ex}\big(\widetilde{R}_1m(R_2.)\big)=E_{ex}(m)$ where $\widetilde{R}_2$ \index{$\widetilde{R}$}is the $3\times 3$-target rotation matrix leaving invariant the component along $e_3$, that is
\begin{equation*}
    \widetilde{R}_2 = \begin{pmatrix}
        R_2 & 0 \\
        0 & 1
    \end{pmatrix}.
\end{equation*}

However, this is not the case for the term $E_{DMI}$, which accounts for an important spin-coupling effect, originating from the chiral nature of the interactions in magnets. The Dzyaloshinskii-Moriya term $E_{DMI}$ is nevertheless invariant under the equivariant action of $SO(2)$: namely, for any rotation matrix $R\in SO(2)$, we have 
\begin{equation*}
    E_{DMI}\big(\widetilde{R}^{-1}m(R \cdot)\big)=E_{DMI}(m).
\end{equation*}

Based on this observation, we define chiral axisymmetric maps of degree $N \in \m Z$ as maps $m: \m R^2 \to \m S^2$ which are $SO(2)$-equivariant  of degree $N$, that is, for any $R\in SO(2)$ and $x\in\R^2$,
\begin{equation}\label{symétrie de degré N}
    m(Rx)=\widetilde{R}^{-N}m(x)= \begin{pmatrix}
        R^{-N} & 0 \\
        0 & 1
    \end{pmatrix}m(x),
\end{equation}
and also satisfy:
\begin{equation}\label{combined reflection 2}
    m(Fx)=\widetilde{-F}m(x) \quad \text{with } F=\begin{pmatrix}
    1 & 0\\
    0 & -1
\end{pmatrix}. \index{$F$}
\end{equation}

The terminology \textit{chiral axisymmetric skyrmions} is due to the chiral properties of the energy $E_{DMI}$.

We claim that chiral axisymmetric maps take the form
\begin{equation}\label{définition skyrmion_01}
    m^{(N,\gamma)}_*(r,\psi)=\begin{pmatrix}
        \sin\big(\Theta(r)\big)\cos(\gamma -N\psi)\\
        \sin\big(\Theta(r)\big)\sin(\gamma -N\psi)\\
        \cos\big(\Theta(r)\big)
    \end{pmatrix} ,
\end{equation}
with helicity\index{helicity} $\gamma\in\frac{\pi}{2}+\R/\pi\Z$, and a with profile\footnote{This profile satisfies the ODE in Proposition~\ref{proposition: propriétés du skyrmion} if $m^{(N,\gamma)}_*(r,\psi)$ is a critical point of the energy.} $\Theta$.

Indeed, let $m^{(N,\gamma)}_*=(m_1,m_2,m_3):\R^2\rightarrow \S^2$ satisfy~\eqref{symétrie de degré N}. We base our argument on the proof of Lemma~A.2 in~\cite{IgnatNguyen1}. Arguing as there, we conclude that $m_3$ is radially symmetric: let $g_1:\R_+\rightarrow\R$ be a function such that $m_3(x)=g_1(|x|)$. We then define for any $r\geq 0$, the function $\Theta(r)=\arccos\big(g_1(r)\big)$. We now consider $m_1$ and $m_2$. Let
\begin{equation*}
    R_\psi = \begin{pmatrix}
        \cos\psi & -\sin\psi\\
        \sin\psi & \cos\psi
    \end{pmatrix},
\end{equation*}
denote the rotation matrix of angle $\psi$. Writing $(e_1,e_2)$ the canonical basis of $\R^2$, we define $\big(a(r),b(r)\big):=\big(m_1(re_2),m_2(re_2)\big)$. Then, denoting $x=(r\cos\psi,r\sin\psi)\in \R^2$ and  applying~\eqref{symétrie de degré N} with $R=R_{\psi-\frac{\pi}{2}}$ implies that for any $(r,\psi),\ Rre_2=x$ and
\begin{equation*}
    m^{(N,\gamma)}_*(r,\psi)=\begin{pmatrix}
        a(r)\cos \big(N\psi-\frac{N\pi}{2}\big) +b(r)\sin  \big(N\psi-\frac{N\pi}{2}\big) \\
        -a(r)\sin  \big(N\psi-\frac{N\pi}{2}\big)  + b(r)\cos  \big(N\psi-\frac{N\pi}{2}\big) \\
        \cos\Theta(r)
    \end{pmatrix}.
\end{equation*}

The symmetry~\eqref{combined reflection 2} implies, for any $(r,\psi)$, the identity
\begin{align*}
    &\begin{pmatrix}
        a(r)\cos  \big(N\psi+\frac{N\pi}{2}\big) -b(r)\sin  \big(N\psi+\frac{N\pi}{2}\big) \\
        a(r)\sin  \big(N\psi+\frac{N\pi}{2}\big) + b(r)\cos  \big(N\psi+\frac{N\pi}{2}\big)\\
        \cos\Theta(r)
    \end{pmatrix}\\
    &\quad\quad \quad =\begin{pmatrix}
        -a(r)\cos  \big(N\psi-\frac{N\pi}{2}\big) -b(r)\sin \big(N\psi-\frac{N\pi}{2}\big) \\
        -a(r)\sin \big(N\psi-\frac{N\pi}{2}\big) + b(r)\cos\big(N\psi-\frac{N\pi}{2}\big) \\
        \cos\Theta(r).
    \end{pmatrix}
\end{align*}

If $N$ is odd (resp. $N$ is even), this implies that $b\equiv 0$ (resp. $a\equiv 0$). At last, using that $m^{(N,\gamma)}_*(r,\psi)$ lies on the sphere leads to $a(r)=\sin\Theta(r)$ (resp. $b(r)=\sin\Theta(r)$) and we obtain the desired form with $\gamma = \frac{N\pi}{2}$ (resp. $\gamma =  \frac{\pi}{2}+ \frac{N\pi}{2}$).

\begin{rem}\label{remarque: Néel skyrmion}
The case corresponding to a helicity\index{helicity} $\gamma\in\R/\pi\Z$ and $N=-1$ in~\eqref{définition skyrmion_01} is $SO(2)$-equivariant in the sense given previously, but the other symmetry now write $m(Fx)=\widetilde{F}m(x)$ instead of \eqref{combined reflection 2}. For such $\gamma$, the mapping $m^{(-1,\gamma)}_*$ associated with the profile $\theta$ provided by Proposition~\ref{proposition: propriétés du skyrmion} is also a critical point, and is called the Néel skyrmion. In either case, a straightforward computation shows that, for $N=-1$,
\begin{equation*}
    E_{DMI}\big(m_*^{(N,\gamma)}\big)=-\frac{4\pi}{3}\sin(\gamma)\int_0^{\infty}\sin^3\big(\theta(r)\big)\,dr.
\end{equation*}
Given the shape of $\theta$, the integral in this identity is positive, so that $ E_{DMI}\big(m_*^{(-1,\gamma)}\big)$ is strictly minimized at $\gamma=\frac{\pi}{2}$. In particular, the Néel skyrmion, corresponding to $\gamma\in\{0,\pi\}$, can never be a minimizer of $E_{DMI}$ among configurations of degree $N=-1$. Rather, it sits at the opposite extremum from the Bloch skyrmion\index{Bloch/Néel skyrmion}. This suggests that the Néel skyrmion is an unstable critical point of the full energy rather than a stable one.
\end{rem}

\section{Estimates for high order terms}\label{appendix: estimates for high order terms}

Due to the lack of Sobolev embeddings in the two dimensional framework, one cannot hope for proper estimates at the $H^1$ level. However, by combining appropriate Gagliardo-Nirenberg type estimates and the Sobolev embeddings $H^1(\R^2)\hookrightarrow L^p(\R^2)$ for any $p\in [2,+\ii)$ and $H^2(\R^2)\hookrightarrow L^\ii(\R^2)$, we shall obtain a sufficient control to prove the crucial estimates in the previous subsections. These controls comprise the following lemmas, giving estimates of integrated quantities that behave as $\mathcal{O}_k^l(\eta)$\index{$\mathcal{O}_k^l(\eta)$} in terms of the Sobolev norms of $\eta$.

\begin{lem}\label{lemme: controle int O^l_1 et O^l_2}
Let $\eta\in H^2(\R^2)$ and an integer $l\geq 2$. We have\index{$\lesssim$}\index{$\mathcal{O}_k^l(\eta)$}
\begin{equation*}
    \left|\int_{\R^2}\mathcal{O}_1^l(\eta)\right|\lesssim \norm{\eta}{H^1}^l,
\end{equation*}

and
\begin{equation*}
    \left|\int_{\R^2}\mathcal{O}_2^l(\eta)\right|\lesssim \norm{\eta}{H^1}^2 \norm{\eta}{L^\ii}^{l-2} .
\end{equation*}

\end{lem}

\begin{proof}
    The proof stems directly from the Sobolev embedding $H^1(\R^2)\hookrightarrow L^p(\R^2)$ for any finite $p\geq 2$.
\end{proof}

\begin{rem}\label{remarque: Gagliardo Nirenberg}
Possibly using the Gagliardo-Nirenberg inequality $\normLii{\eta}\lesssim \norm{\eta}{H^1}^{\frac{1}{2}}\norm{\eta}{H^2}^{\frac{1}{2}}$, the above lemma could be refined by reducing the exponent of the higher Sobolev norm. It is not necessary to do so in the above lemma, but we use these refinements in the next one.
\end{rem}

\begin{lem}\label{lemme: controle int O^l_3 et O^l_4}
Let $\eta\in H^3(\R^2)$. We have\index{$\mathcal{O}_k^l(\eta)$}\index{$\lesssim$}
    \begin{equation*}
        \left|\int_{\R^2}\mathcal{O}_4^2(\eta)\right|\lesssim \norm{\eta}{H^2}^2,
    \end{equation*}
and
    \begin{equation*}
        \left|\int_{\R^2}\mathcal{O}_5^2(\eta)\right|\lesssim \norm{\eta}{H^2}\norm{\eta}{\dot{H}^3}.
    \end{equation*}

Using the Gagliardo-Nirenberg inequality, we obtain
\begin{equation*}
    \left|\int_{\R^2}\mathcal{O}_5^3(\eta)\right|\lesssim \big(\normLii{\eta}+\normHun{\eta}\big)\norm{\eta}{H^2}\norm{\eta}{\dot{H}^3},
\end{equation*}
    \begin{equation*}
        \left|\int_{\R^2}\mathcal{O}_6^3(\eta)\right|\lesssim \norm{\eta}{L^\ii}\norm{\eta}{\dot{H}^3}^2+ \norm{\eta}{H^1}^{\frac{1}{2}}\norm{\eta}{H^2}\norm{\eta}{\dot{H}^3}^{\frac{3}{2}} + \norm{\eta}{H^2}^2 \norm{\eta}{\dot{H}^3}, 
    \end{equation*}
    as well as
    \begin{align*}
        \left|\int_{\R^2}\mathcal{O}_6^4(\eta)\right|\lesssim \normLii{\eta}^2\norm{\eta}{\dot{H}^3}^2 +& \normLii{\eta}\normHun{\eta}^\frac{1}{2} \norm{\eta}{H^2}\norm{\eta}{\dot{H}^3}^\frac{3}{2} \\
        &+ \big(\normLii{\eta}+\normHun{\eta}\big)\norm{\eta}{H^2}^2\norm{\eta}{\dot{H}^3} . 
    \end{align*}
\end{lem}

\begin{proof}
We only deal with the last estimate, the proofs of the others can be derived by adapting the following estimates. After appropriate integrations by parts, the potential terms appearing in $\mathcal{O}_6^4(\eta)$ are of the form $(\partial_x^3\eta)^2\eta^2,\partial_x^3\eta\partial_x^2\eta\partial_x\eta\eta$ and $\partial_x^3\eta(\partial_x\eta)^3$. The first inequality
\begin{equation*}
    \left|\int_{\R^2}(\partial_x^3\eta)^2\eta^2\right|\lesssim \normLii{\eta}^2\norm{\eta}{\dot{H}^3}^2,
\end{equation*}
is straightforward. As for the next term, we write
\begin{equation*}
    \left|\int_{\R^2}\partial_x^3\eta\partial_x^2\eta\partial_x\eta\eta\right|\leq \norm{\eta}{\dot{H}^3}\norm{\partial_x^2\eta}{L^4}\norm{\partial_x\eta}{L^4}\normLii{\eta},
\end{equation*}
and we conclude by invoking the two-dimensional Gagliardo-Nirenberg inequality, $\norm{\eta}{L^4}\lesssim \norm{\partial_x\eta}{L^2}^{\frac{1}{2}}\norm{\partial_x\eta}{L^2}^{\frac{1}{2}}$. We conclude by using the same Gagliardo-Nirenberg inequality in the estimate
\begin{align*}
    \left|\int_{\R^2}\partial_x^3\eta (\partial_x\eta)^2(\partial_x\eta+\eta)\right|&\leq \left|\int_{\R^2}(\partial_x^2\eta )^2(\partial_x\eta)^2\right| \leq \norm{\partial_x^2\eta}{L^4}^2\norm{\partial_x\eta}{L^4}^2.
\end{align*}
\end{proof}

Moreover, the properties of $f$ can be cleverly combined so that we deduce the following result.
\begin{prop}\label{proposition: propriétés de f}
Up to increasing the value of $B_*$, we have:
\begin{enumerate} 
\item for any $x\in\R_+^*$,   $x\leq B_*f(x,x)$.
\item for any $x,y,z,t\in\R_+$, we have 
\begin{equation*}
    f(x+y,z+t)\leq B_*\big( f(x,y) + f(x,z)+f(y,z)+f(y,t)\big).
\end{equation*}
\item for any $\eta\in H^2(\R^2)$,
    \begin{equation*}
        \normLii{\eta}+\normHun{\eta}\leq  B_*f(\normHun{\eta},\norm{\eta}{H^2}) .
    \end{equation*}
\end{enumerate}
\end{prop}

\begin{proof}
The first property can be derived by taking $x=y=\frac{1}{\lambda}>0$ in~\eqref{2eme hypothese f} for any $x>0$. The second is deduced using~\eqref{1ere hypothese f} with the fact that $y+z\leq 2\max(y,z)$.

The last one is a combination of the properties of the previous both and from~\eqref{3eme hypothese f}.
\end{proof}

\section{Proofs of the local well-posedness result}\label{appendix: LWP}

Our goal here is to prove Theorems \ref{théorème: local well-posedness} and \ref{théorème: LWP of the Cauchy problem}. We start with the former. In the sequel, we shall assume that $d\in\{2,3\}$ and $s\geq 2$. 

Let us fix\index{$\mathcal{Z}^{s}(\R^d)$} $\mathfrak{m}\in \mathcal{Z}^{s}(\R^d)$ such that $\mathfrak{m}\times \heff(\mathfrak{m})\in H^{s-1}(\R^d)$. We write $Lw=-\Delta w+Nw = -\Delta w +2\grad\times w + hw$ and $N_* \mathfrak{m}:= 2 \grad\times \mathfrak{m} + h(\mathfrak{m}-e_3)$ so that $N(\mathfrak{m}+w)=N_*\mathfrak{m} + Nw$.

In particular $m=\mathfrak{m}+w$ solves~\eqref{llg} if and only if (formally) $w$ solves
\begin{equation}\label{llg on w}
    \partial_tw=\alpha\Delta w + (\mathfrak{m}+w)\times \Delta w + Z(w),
\end{equation}
with $Z(w)=Z_0 + Z_1(w) + \alpha Z_2(w)$, and where each term in the previous sum is defined in the following:
\begin{equation}\label{définition Z_0}
    Z_0=\mathfrak{m}\times \heff(\mathfrak{m}) + \alpha\Big(\big(|\grad \mathfrak{m}|^2  +\mathfrak{m}.N_*(\mathfrak{m})\big)\mathfrak{m} +\heff(\mathfrak{m})\Big) ,
\end{equation}

\begin{equation*}
    Z_1(w)=w\times \heff(\mathfrak{m}) - (\mathfrak{m}+w)\times Nw, 
\end{equation*}
and
\begin{align*}
    Z_2(w)=\Big(2\grad \mathfrak{m}:\grad w +&|\grad w|^2 + w.N_*(\mathfrak{m}+w) + \mathfrak{m}.Nw\Big) (\mathfrak{m}+w)\\
    &+|\grad \mathfrak{m}|^2 w+\mathfrak{m}.N_*\mathfrak{m}\ w - Nw.
\end{align*}

An important fact is the following claim.
\begin{claim}\label{claim: Z_0 dans H^s-1 si et seulement si mathfrak m in Z^s et mathfrak m..}
Under the assumptions in Theorem~\ref{théorème: LWP of the Cauchy problem}, $Z_0\in H^{s-1}(\R^d)$.
\end{claim} 

\begin{proof}
By definition $Z_0 = \mathfrak{m}\times\heff(\mathfrak{m}) - \alpha \mathfrak{m}\times \big(\mathfrak{m}\times\heff(\mathfrak{m})\big)$ and~\eqref{définition Z_0} is obtained using $|\mathfrak{m}|^2=1$ and differentiating two times this identity.\\
Now take $\gamma\in\N^d$ with $|\gamma| = s-1$, we have by hypothesis, that for any $\widetilde{\gamma}\leq\gamma,\partial^{\widetilde{\gamma}}\big(\mathfrak{m}\times \heff (\mathfrak{m})\big)\in L^2(\R^d)$. Let us check that it is also the case for $\partial^{\widetilde{\gamma}}\Big(\mathfrak{m}\times\big(\mathfrak{m}\times \heff (\mathfrak{m})\big)\Big)$. Regarding the $0$-order term, we have 
\[ \normLdeux{\mathfrak{m}\times\big(\mathfrak{m}\times \heff (\mathfrak{m})\big)}\leq \normLii{\mathfrak{m}}\normLdeux{\mathfrak{m}\times \heff (\mathfrak{m})}\leq \norm{\mathfrak{m}\times \heff (\mathfrak{m})}{H^{s-1}}. \]
 On the other hand, we estimate
\begin{align*}
    &\normLdeux{\partial^{\gamma}\Big(\mathfrak{m}\times\big(\mathfrak{m}\times \heff (\mathfrak{m})\big)\Big)}\\
    &\lesssim \normLdeux{\partial^\gamma \mathfrak{m} \times \big(\mathfrak{m}\times \heff(\mathfrak{m})\big)}+ \sum_{\substack{\gamma^{(1)}+\gamma^{(2)}=\gamma \\ \gamma^{(2)}\neq 0}}\normLii{\partial^{\gamma^{(1)}}\mathfrak{m}}\normLdeux{\partial^{\gamma^{(2)}}\big(\mathfrak{m}\times \heff(\mathfrak{m})\big)}\\
    &\lesssim \norm{\partial^\gamma \mathfrak{m} }{L^4}\norm{\mathfrak{m}\times \heff(\mathfrak{m})}{L^4} + \sum_{\substack{\gamma^{(1)}+\gamma^{(2)}=\gamma \\ \gamma^{(2)}\neq 0}}\norm{\partial^{\gamma^{(1)}}\mathfrak{m}}{H^2}\norm{\mathfrak{m}\times \heff(\mathfrak{m})}{H^{s-1}}\\
    &\lesssim \norm{\partial^\gamma \mathfrak{m} }{H^1}\norm{\mathfrak{m}\times \heff(\mathfrak{m})}{H^1}+  \norm{\grad \mathfrak{m}}{H^{s-1}}\norm{\mathfrak{m}\times \heff(\mathfrak{m})}{H^{s-1}} \\
    &\lesssim \norm{\mathfrak{m}}{\mathcal{Z}^s}\norm{\mathfrak{m}\times \heff(\mathfrak{m})}{H^{s-1}}.
\end{align*}
\end{proof}

We use a standard regularization method to solve~\eqref{llg on w}, that is, we consider the initial value problem
\begin{equation}\label{llg sur w regularisation}
    \left\{\begin{array}{ll}
    \partial_t w = -\beta \Delta^2 w +\alpha\Delta w +(\mathfrak{m}+w)\times \Delta w + Z(w), \\ 
    w(0)=w_0 \in H^s(\R^d).
    \end{array}
    \right.
\end{equation}

Regarding the previous equation, we can proceed to a standard fixed point argument using the Duhamel formula. Writing $\widetilde{Z}(w)=\alpha \Delta w + (\mathfrak{m}+w)\times \Delta w +  Z(w)$, we consider the mapping 
\begin{equation}
    \Phi(w):=e^{-t\beta \Delta^2}w_0 + \int_0^t e^{-(t-\tau)\beta \Delta^2}\widetilde{Z}\big(w(\tau)\big)d\tau.
\end{equation} 

Then, we adapt the argument in~\cite{shimizu2026} (see Appendix~A), we only recall the main steps of this argument and provide the necessary details so that the present argument is self-contained.

The main discrepancy from the latter proof is the source term $Z_0$. This term does not influence the contractive property of $\Phi$, which remains true. Regarding the target set of $\Phi$, we invoke both very general estimates $(A.5)$ and $(A.6)$ in~\cite{shimizu2026}, according to which there exists a constant $C_\beta$ and for any $k\in [0,2]$, a constant $C_{\beta,k}$ such that for any $t>0$ and $f\in L^2(\R^d) $, we have
\begin{equation}\label{premiere inégalité}
    \norm{e^{-t\beta\Delta^2}f}{L^\ii_t L^2_x} + \norm{\Delta e^{-t\beta\Delta^2}f}{L^2_{t,x}} \leq C_\beta \norm{f}{L^2_x},
\end{equation}
and
\begin{equation}\label{deuxieme inegalité}
    \norm{e^{-t\beta\Delta^2}\Pi _H f}{L^\ii_t L^2_x}\leq C_{\beta,k}t^{-\frac{k}{4}}\norm{|\grad|^{-k}\Pi_H f}{L^2_x}.
\end{equation}

In our case, considering
\begin{equation*}
    G(t,x)=\int_0^t e^{-(t-\tau)\beta\Delta^2}\Pi_HF(\tau)d\tau,
\end{equation*}
which is a solution of the equation $\partial_t G + \beta\Delta^2G = \Pi_HF$ and writing $\widetilde{G}(\tau,\xi)$ the space-time Fourier transform, the evolution equation upon $G$ reads $i\widetilde{G} + \beta |\xi|^4\widetilde{G}=\widetilde{\Pi_H F}$, leading to $$\widetilde{G}(\tau,\xi)=\frac{\widetilde{\Pi_H F}}{i\tau + \beta|\xi|^4}.$$ 
Then for $k\in (0,2]$, we infer by the one-dimensional Sobolev embedding $\dot{H}_t^{\frac{2-k}{4}}\hookrightarrow L^{\frac{4}{k}}_t$ that \begin{equation*}
    \norm{\Delta G}{L^{\frac{4}{k}}_tL^2_x}\lesssim\norm{\partial_t^{\frac{2-k}{4}}\Delta G}{L^2_tL_x^2} \lesssim\norm{\dfrac{\tau^{\frac{2-k}{4}}|\xi|^{2+k}}{i\tau+\beta|\xi|^4} |\xi|^{-k}\widetilde{\Pi_H F}}{L^2_\tau L^2_\xi}.
\end{equation*}

The standard convexity inequality $\tau^{\frac{2-k}{4}}|\xi|^{2+k}\leq \frac{(2-k)\tau}{4} + |\xi|^4$ provides
\begin{equation}\label{troisieme inegalité}
    \norm{\Delta G}{L^{\frac{4}{k}}_tL^2_x}\lesssim \norm{|\grad|^{-k}\Pi_HF}{L^2_tL_x^2}.
\end{equation}

By separating in high and low frequence, and by using all three previous estimates~\eqref{premiere inégalité},~\eqref{deuxieme inegalité} and~\eqref{troisieme inegalité}, there exists a constant $D_\beta>0$ such that for any $t\in [0,T]$, 
\begin{align*}
    \norm{\int_0^t e^{-(t-\tau)\beta \Delta^2} Z_0 d\tau}{H^s_x}&\lesssim C_\beta \int_0^T\norm{\Pi_L Z_0}{L^2_x} + C_\beta\int_0^t (t-\tau)^{-\frac{1}{2}}\norm{\Pi_H Z_0}{H_x^{s-2}}\\
    &\lesssim D_\beta (T^{\frac{1}{2}}+T)\norm{Z_0}{L^\ii_TH^{s-2}_x},
\end{align*}
and similarly,
\begin{equation*}
    \norm{\Delta\int_0^t e^{-(t-\tau)\beta \Delta^2} Z_0 d\tau}{H^s_x}\lesssim D_\beta (T^{\frac{1}{2}}+T^{\frac{3}{2}})\norm{Z_0}{L^\ii_TH^{s-2}_x}.
\end{equation*}

Both previous estimates eventually provide two constants $K\geq 1$ and $C>0$, which neither depend on $\gm$, nor $T$, such that 
\begin{align*}
\MoveEqLeft  \norm{\int_0^t e^{-(t-\tau)\beta \Delta^2} Z_0 d\tau}{L^\ii_T H^s_x}+\norm{\Delta\int_0^t e^{-(t-\tau)\beta \Delta^2} Z_0 d\tau}{L^2_T H^s_x} \\
    & \le  C(T+T^\frac{1}{2}+T^\frac{3}{2})\norm{Z_0}{L^\ii_T H^{s-2}_x}
    \le C(T^\frac{1}{2}+T^\frac{3}{2})\big(\norm{\mathfrak{m}}{{\mathcal{Z}^{s}}}+\norm{\mathfrak{m}}{\mathcal{Z}^{s}}^K\big),
\end{align*}

hence, for possibly larger constants $C$ and $K$,
\begin{align*}
    \norm{\Phi(w)}{L^\ii_T H^s_x}+\norm{\Delta\Phi(w)}{L^2_T H^s_x}\leq C(T^\frac{1}{2}+T^\frac{3}{2})\big(\norm{\mathfrak{m}}{{\mathcal{Z}^{s}}}+\norm{\mathfrak{m}}{\mathcal{Z}^{s}}^K\big).
\end{align*}

In a nutshell, by setting $T=(1+\norm{\gm}{\mathcal{Z}^s})^{\frac{2(K-1)}{3}}$ and $M=C(1+\norm{\gm}{\mathcal{Z}^s})$, we deduce that for $\mathfrak{m}\in \mathcal{Z}^s(\R^d)$, the mapping $\Phi$ is a contraction on 
\begin{align*}
    &X_T:=\\
    \Big \{ & w\in \q C^0\big([0,T],H^s(\R^d)\big)\cap L^2\big([0,T],H^{s+2}(\R^d)\big)\Big| \norm{w}{L^\ii_T H^s_x} + \norm{\Delta w}{L^2_T H^s_x}\leq M\Big\}.
\end{align*}

Thus, for any $\beta >0$, there exist $T_{\max}^\beta>0$ and a solution 
\[ w^\beta\in \q C^0\big((0,T^\beta_{\max}),H^s(\R^d)\big)\cap L^2\big((0,T^\beta_{\max}),H^{s+2}(\R^d)\big) \] to~\eqref{llg sur w regularisation} with initial data $w_0$. Now, we find uniform estimates with respect to $\beta$, which will allow us to pass to the limit $\beta\rightarrow 0$ and obtain a solution to~\eqref{llg on w}. We first state the appropriate uniform estimates associated with the fourth-order regularization~\eqref{llg sur w regularisation}.

\begin{lem}\label{lemme: normw^betaL^ii_T H^s_x +normgrad w^betaL^2_TH^s_x+ betanormDelta w^betaL^2_TH^s_xleq C}
    Let $w_0\in H^s(\R^d)$ and $w^\beta\in \q C^0\big([0,T^\beta_{\max}),H^s(\R^d)\big)$ the unique solution to~\eqref{llg on w}. There exist positive constants $C$ and $T$, only depending on $w_0$ such that any $\beta >0,0<T< T^\beta_{\max}$ and
    \begin{equation}
        \norm{w^\beta}{L^\ii_T H^s_x} + \norm{\grad w^\beta}{L^2_TH^s_x}+ \beta\norm{\Delta w^\beta}{L^2_TH^s_x}\leq C.
    \end{equation}
\end{lem}

\begin{proof}
For any spatial multi-index $\gamma\in\N^d$ with $|\gamma|\leq s$, we compute
\begin{align}
    \dfrac{1}{2}\dfrac{d}{dt}&\normLdeux{\partial^\gamma w(t)}^2 = -\beta\normLdeux{\Delta \partial^\gamma w(t) }^2 -\alpha\normLdeux{\grad\partial^\gamma w(t)}^2 \label{d/dt norm Ldeux partial gamma w(t)}\\
    & + \psLdeux{(\mathfrak{m}+w)\times \partial^\gamma \Delta w}{\partial^\gamma w}+\psLdeux{\partial^\gamma Z_0}{\partial^\gamma w(t)} + \psLdeux{\widetilde{Z}_{\gamma}(w(t))}{\partial^\gamma w(t)}\notag, 
\end{align}

with \begin{align*}
    \widetilde{Z}_\gamma (w)=\sum_{\substack{\gamma^{(1)}+\gamma^{(2)}=\gamma \\ \gamma^{(1)}\neq 0}}\big(\partial^{\gamma^{(1)}}(\mathfrak{m}+w)\times \Delta \partial^{\gamma^{(2)}}w \big) +\partial^\gamma Z(w).
\end{align*}

\begin{claim}\label{claim: estimée unforme en beta}
There exist $0<T<T_{\max}^\beta,K\geq 2$ and $C>0$, that do not depend on $\beta$, such that     \begin{align*}
    \dfrac{1}{2} \norm{\partial^\gamma w}{L^\ii_T L^2_x}^2 &+  \alpha\norm{\grad\partial^\gamma w}{L^2_TL^2_x}^2 + \beta \norm{\Delta \partial^\gamma w }{L^2_TL^2_x}^2 \\
    &\leq \dfrac{1}{2}\norm{\partial^\gamma w(0)}{ L^2_x}^2 + C\varepsilon(T)\bigg(\norm{w}{L^\ii_TH^s_x} +\norm{w}{L^\ii_TH^s_x}^K \\
    & \quad \quad \quad\quad \quad \quad + \norm{\grad w}{L^2_TH^s_x}^2 +\int_0^T\psLdeux{\partial^\gamma Z_0}{\partial^\gamma w(t)}dt\bigg),
\end{align*}
with $\varepsilon(T)\rightarrow 0$ as $T\rightarrow 0$.

\end{claim}

\begin{proof}[Proof of Claim~\ref{claim: estimée unforme en beta}]

The same cancelation than in~\cite{shimizu2026} occurs for the term involving $(\mathfrak{m}+w)\times \Delta\partial^\gamma w$ and we similarly obtain
\begin{equation*}
    \psLdeux{(\mathfrak{m}+w)\times \Delta\partial^\gamma w}{\partial^\gamma w}=-\sum_{j=1}^d\psLdeux{\partial_j(\mathfrak{m}+w) \times \partial_j\partial^\gamma w}{\partial^\gamma w}.
\end{equation*}

Thus, integrating~\eqref{d/dt norm Ldeux partial gamma w(t)} between $0$ and $T>0 $ leads to
\begin{align}
\MoveEqLeft[4]    \dfrac{1}{2}\normLdeux{\partial^\gamma w(T)}^2 -\dfrac{1}{2} \normLdeux{\partial^\gamma  w(0)}^2  +\int_0^T\Big(\beta\normLdeux{\Delta \partial^\gamma w(t) }^2 + \alpha\normLdeux{\grad\partial^\gamma w(t)}^2\Big)dt \nonumber \\
    &=-\sum_{j=1}^d\int_0^T\psLdeux{\partial_j(\mathfrak{m}+w) \times \partial_j\partial^\gamma w}{\partial^\gamma w}\label{terme un dans l'intégration}\\
    &\quad +\int_0^T\psLdeux{\partial^\gamma Z_0}{\partial^\gamma w(t)}dt +\int_0^T\psLdeux{\partial^\gamma Z(w)}{\partial^\gamma w}dt \label{terme deux dans l'intégration}\\
    &\quad +\sum_{\substack{\gamma^{(1)}+\gamma^{(2)}=\gamma\\ \gamma^{(1)}\neq 0}}\int_0^T \psLdeux{\big(\partial^{\gamma^{(1)}}(\mathfrak{m}+w)\times \Delta \partial^{\gamma^{(2)}}w \big)}{\partial^\gamma w(t)}\label{terme trois dans l'intégration}.
\end{align}

It remains to control~\eqref{terme un dans l'intégration},~\eqref{terme deux dans l'intégration} and~\eqref{terme trois dans l'intégration} as desired. Regarding~\eqref{terme un dans l'intégration}, we write for any $j$ that 
\begin{align*}
    &\bigg|\int_0^T\psLdeux{\partial_j(\mathfrak{m}+w) \times \partial_j\partial^\gamma w}{\partial^\gamma w}\bigg|\\
    &\leq \normLii{\grad \mathfrak{m}}\int_0^T \norm{\grad w}{H^s_x}\norm{w}{H^s_x} + \int_0^T\norm{\grad w}{H^s_x}\norm{\grad w}{L^4_x}\norm{\partial^\gamma w}{L^4_x}\\
    &\leq \normLii{\grad \mathfrak{m}}T^\frac{1}{2}\norm{w}{L^\ii_TH^s_x}\norm{\grad w}{L^2_T H^s_x}+ \int_0^T\norm{\grad w}{H^s_x}\norm{ w}{H^s_x}\norm{w}{H^s_x}^{1-\frac{d}{4}}\norm{\grad w}{H^s_x}^\frac{d}{4}\\
    &\leq \normLii{\grad \mathfrak{m}}T^\frac{1}{2}\Big(\norm{w}{L^\ii_TH^s_x}^2+\norm{\grad w}{L^2_T H^s_x}^2\Big) + T^{\frac{1}{2}-\frac{d}{8}}\norm{w}{L^\ii_TH^s_x}^{2-\frac{d}{4}}\norm{\grad w}{L^2_T H^s_x}^{1+\frac{d}{4}}\\
    &\leq C_1(\mathfrak{m})\big(T^\frac{1}{2} + T^{\frac{1}{2}-\frac{d}{8}}\big)\Big(\norm{w}{L^\ii_TH^s_x}^2+\norm{\grad w}{L^2_T H^s_x}^2+\norm{\grad w}{L^2_T H^s_x}^K\Big),
\end{align*}
avec $K\geq 2$. As for~\eqref{terme deux dans l'intégration}, it is sufficient to observe that 
\begin{align*}
    \int_0^T\psLdeux{\partial^\gamma Z_0}{\partial^\gamma w(t)}\leq \norm{\partial^\gamma Z_0}{H^{-1}_x}\int_0^T\norm{\grad w}{H^s_x}\leq \norm{\mathfrak{m}}{\mathcal{Z}^{s+1}_x}T^\frac{1}{2}\norm{\grad w}{L^2_TH^s_x}.
\end{align*}

The last term~\eqref{terme trois dans l'intégration} can be dealt with exactly as in~\cite{shimizu2026}.
\end{proof}

Then, since the estimate in Claim~\ref{claim: estimée unforme en beta} is uniform in $\beta$, Lemma~\ref{lemme: normw^betaL^ii_T H^s_x +normgrad w^betaL^2_TH^s_x+ betanormDelta w^betaL^2_TH^s_xleq C} follows directly by implementing a continuation argument.
\end{proof}

Then, 
we can extract a subsequence $(\beta_n)$ such that $w^{\beta_n}$ converges in a certain sense towards $w\in \q C^0\big([0,T_{\max},H^s(\R^d)\big)\cap L^2\big([0,T_{\max}),H^{s+1}(\R^d)\big)$. This concludes the existence part in the local-wellposedness result for~\eqref{llg on w}. In addition, taking to solution $w^1$ and $w^2$ two solutions associated with initial data $w_0^1$ and $w_0^2$ in $H^s(\R^d)$, then the source term $Z_0$ does not intervene in the equation satisfied by the evolution equation satisfied by $w^1-w^2$ and 
one can find $\widetilde{T}< T$ small enough such that 
\begin{equation*}
    \norm{ w^1-w^2}{L^\ii_{\widetilde{T}}H^{s-1}_x}\leq C\norm{w_0^1-w_0^2}{H^{s-1}}, 
\end{equation*}

where $C$ depends on $\norm{w^j}{L^\ii_{\widetilde{T}} H^{s}_x}$ and $\norm{w^j}{L^2_{\widetilde{T} }H^{s+1}_x}$ for $j\in\{1,2\}$, which simultaneously implies the uniqueness of the solution of~\eqref{llg on w} with prescribed initial data $w_0\in H^s(\R^d)$ on a maximal interval $[0,T_{\max})$, and the continuity of the flow with respect to $H^{s-1}$.

We can now prove the main local-wellposedness result.
\begin{proof}[Proof of Theorem~\ref{théorème: local well-posedness}]

Take $\mathfrak{m}\in\mathcal{Z}^s(\R^d)$ satisfying~\eqref{Z_0 in H^s-1} and an initial data $w_0\in H^s(\R^d)$. The uniqueness is a straightforward consequence of the uniqueness in the Cauchy problem~\ref{llg on w}.

We prove that $ w\in \q C^1\big([0,T_{\max}), H^{s-2}(\R^d)\big)$. To do this, and by~\eqref{llg on w}, it is sufficient to check that $\Delta w, (\mathfrak{m}+w)\times \Delta w $ and $Z(w)$ lie in $\q C^0\big([0,T_{\max}), H^{s-2}(\R^d)\big)$. We shall use the previous construction, in particular the fact that $w\in \q C^0\big([0,T_{\max}), H^{s}(\R^d)\big)$. We deal first with the quasilinear term. Let $\gamma\in\N^d$ satisfying $|\gamma| = s-2$. We have 
\begin{align*}
    \norm{(\mathfrak{m}+w)\times \Delta w}{ H^{s-2}_x}&\lesssim \sum_{\gamma^{(1)}+\gamma^{(2)}=\gamma}\norm{\partial^{\gamma^{(1)}}\Delta w}{L^2_x}\norm{\partial^{\gamma^{(2)}}(\mathfrak{m}+ w)}{L^\ii_x}\\
    & \lesssim \sum_{\gamma^{(1)}+\gamma^{(2)}=\gamma}\norm{\partial^{\gamma^{(1)}}w}{H^2_x}\Big(\norm{\partial^{\gamma^{(2)}}\mathfrak{m}}{L^\ii_x}+\norm{\partial^{\gamma^{(2)}}w}{L^\ii_x}\Big)\\
    &\lesssim \sum_{\gamma^{(1)}+\gamma^{(2)}=\gamma}\norm{w}{H^s_x}\Big(\norm{\partial^{\gamma^{(2)}}\mathfrak{m}}{L^\ii_x}+\norm{w}{H^s_x}\Big)
\end{align*}

If $|\gamma^{(2)}|=0$, $\norm{\partial^{\gamma^{(2)}}\mathfrak{m}}{L^\ii_x}\leq \norm{\mathfrak{m}}{L^\ii_x}\leq \norm{\mathfrak{m}}{\mathcal{Z}^s_x}$ and if $|\gamma^{(2)}|>0$, 
\[ \norm{\partial^{\gamma^{(2)}}\mathfrak{m}}{L^\ii_x}\lesssim \norm{\partial^{\gamma^{(2)}}\mathfrak{m}}{H^2_x}\leq\norm{\nabla \mathfrak{m}}{H^{s-1}_x}\leq \norm{\mathfrak{m}}{\mathcal{Z}^s_x}. \]
 Therefore, $(\mathfrak{m}+w)\times \Delta w\in \q C^0\big([0,T_{\max}), H^{s-2}(\R^d)\big)$. Similarly, we can prove that $Z_1(w)$ and $Z_2(w)$ belong to $\q C^0\big([0,T_{\max}), H^{s-2}(\R^d)\big)$. As a matter of example, we deal with one single term in $Z_2(w)$. Namely, take $\gamma$ as previously, we have
\begin{align*}
    \norm{|\grad w|^2 (\mathfrak{m}+w)}{H^{s-2}_x}&\lesssim \sum_{\gamma^{(1)}+\gamma^{(2)}=\gamma}\norm{\partial^{\gamma^{(1)}}(|\grad w|^2)\partial^{\gamma^{(2)}}(\mathfrak{m}+w)}{L^2_x}\\
    &\lesssim \sum_{\gamma^{(1)}+\gamma^{(2)}=\gamma}\norm{\partial^{\gamma^{(1)}}(|\grad w|^2)}{L^2_x}\norm{\partial^{\gamma^{(2)}}(\mathfrak{m}+w)}{L^\ii_x}\\
    &\lesssim \norm{\partial^\gamma\grad w}{L^4_x}^2\Big(\norm{\mathfrak{m}}{\mathcal{Z}^s_x} + \norm{w}{H^s_x}\Big)\\
    &\lesssim \norm{ w}{H^s_x}^2\Big(\norm{\mathfrak{m}}{\mathcal{Z}^s_x} + \norm{w}{H^s_x}\Big).
\end{align*}

As a consequence, $m\in \q C^1\big([0,T_{\max}), H^{s-2}(\R^d)\big)$ and $\partial_t m = \partial_t w$ that belongs to $\q C^0\big([0,T_{\max}), L^2(\R^d)\big)$. We claim that the scalar function $h(t)=\normLdeux{1-|m(t)|^2}^2$ is differentiable. Indeed, since $|\mathfrak{m}|=1$, then $(1-|m(t)|^2)^2\leq 2\big(2|w(t)|^2 + |w(t)|^4\big)\in L^1(\R^d)$, due to the Sobolev embedding $H^1(\R^d)\hookrightarrow L^p(\R^d)$ for any finite $p$. Moreover, \begin{equation*}
    \left|\partial_t \Big[(1-|m(t)|^2)^2\Big]\right| = 4\big(2|w|+|w|^2\big)\big(|\partial_t w(t)|+|\partial_t w(t)||w(t)|\big)\in L^1(\R^d),
\end{equation*}

so that, using~\eqref{llg}, we have
\begin{equation*}
    h'(t)=-4\int_{\R^d}\partial_t m(t).m(t) (1-|m(t)|^2)=0.
\end{equation*}

Since by hypothesis $h(0)=0$, thus $h(t)=0$ for any $t\in [0,T_{\max})$, which concludes the proof of Theorem~\ref{théorème: local well-posedness}.
\end{proof}

\begin{proof}[Proof of Theorem~\ref{théorème: LWP of the Cauchy problem}]
Take $d = 2,\ s\in\{2,3\}$ and $m_0\in H^2_*$. By the properties of the skyrmion, we have $\gm=m_*\in \mathcal{Z}^s(\R^2)$ and satisfies~\eqref{Z_0 in H^s-1}, thus Theorem~\ref{théorème: local well-posedness} applies and implies that there exists a solution that decomposes $m(t)=m_* + w(t)$ with $w\in \q C^0\big([0,T_{\max}),H^s(\R^2)\big)$ and initial data $m(0)=m_0$.

Now, $\partial_t m$ and $\delta E\big(m(t)\big)$ belong to $\q C^0\big([0,T_{\max}),L^2(\R^2)\big)$, thus $t\mapsto E\big(m(t)\big)$ can be differentiated. This justifies the expression in~\eqref{dérivée de E(m(t)) non rigoureuse} and implies 
\begin{equation*}
    \dfrac{d}{dt}\big(E(m(t)\big)=Diss_1\big(m(t)\big).
\end{equation*}

If now the data $m_0$ is in $H^3_*$, then $w\in \q C^0([0,T_{\max}), H^3(\R^2)\big)$, thus $m\in \q C^1([0,T_{\max}), H^1(\R^2)\big)$, and $\delta E\big(m(t)\big) \in \q C^0([0,T_{\max}), H^1(\R^2)\big)$. Therefore, $t\mapsto\normLdeux{\delta E\big(m(t)\big)}$ is differentiable. 
\end{proof}

%
%

\printindex

\bibliographystyle{plain}
\bibliography{bibliography}

\end{document}